\documentstyle[12pt,twoside]{amsart}

\renewcommand{\c}[0]{{\mathbb C}}  

\renewcommand{\o}[0]{{\cal O}} 
\newcommand{\z}[0]{{\mathbb Z}}

\renewcommand{\r}[0]{{\mathbb R}} 

\renewcommand{\a}[0]{{\mathbb A}} 

\newcommand{\h}[0]{{\Bbb H}}
\newcommand{\p}[0]{{\mathbb P}}

\newcommand{\q}[0]{{\mathbb Q}}

\newcommand{\qtq}[1]{\quad\mbox{#1}\quad}
\newcommand{\spec}[0]{\operatorname{Spec}}

\newcommand{\red}[0]{\operatorname{red}}

\newcommand{\proj}[0]{\operatorname{Proj}}

\newcommand{\inter}[0]{\operatorname{Int}}    
\newcommand{\sing}[0]{\operatorname{Sing}}    
\newcommand{\ex}[0]{\operatorname{Ex}}    

\newcommand{\nec}[1]{\overline{NE}({#1})}

\newsymbol\subsetneq 2328
\newsymbol\onto 1310
\def\into{\DOTSB\lhook\joinrel\rightarrow}
\newsymbol\twoheadrightarrow 1310

\newtheorem{thm}{Theorem}[section]

\newtheorem{lem}[thm]{Lemma}
\newtheorem{cor}[thm]{Corollary}

\newtheorem{prop}[thm]{Proposition}

\newtheorem{complement}[thm]{Complement}

\theoremstyle{definition}
\newtheorem{defn}[thm]{Definition}

\newtheorem{say}[thm]{}
\newtheorem{exmp}[thm]{Example}

\newtheorem{rem}[thm]{Remark}          
\newtheorem{ack}{Acknowledgments}        
\newtheorem{notation}[thm]{Notation}

\theoremstyle{remark}

\begin{document}
\bibliographystyle{amsplain}

\title{Real Algebraic Threefolds III.\\ Conic Bundles}
\author{J\'anos Koll\'ar}

\maketitle
\tableofcontents

\section{Introduction}

This paper continues the study of the topology of real algebraic
threefolds begun in \cite{rat1, rat2}, but the current work
is completely independent of the previous ones in its methodology.

The present aim is to understand the topology of a special class of
 real algebraic
threefolds. In many respects, the simplest 3--folds are
of the form $\p^1\times S$ where $S$ is a surface. If one adopts the
birational point of view, it is more natural, and considerably more
general, to investigate those 3-folds 
which admit a map onto a surface  $f:X\to S$ whose general fibers are
smooth rational curves. 
This class of threefolds also appears as one of the 4 possible
outcomes of the minimal model program (cf.\ \cite{KoMo98}). 
From the birational point of view the most interesting 3--folds in this
class are those which map onto a rational surface. 
This is also the natural assumption to make in connection with the Nash
conjecture \cite[p.\ 421]{Nash52}.  Our main theorem gives a nearly complete
description of the possible topological types of the set of real points
of such a threefold. (See (\ref{what.is.rav?}) for 
 terminology and notation.)

\begin{thm}\label{orient.main.thm} Let $X$ be a smooth projective real
algebraic threefold such that the set of real points $X(\r)$ is
orientable.
   Let $f:X\to S$ be a morphism
onto a real algebraic surface  $S$ whose 
 general fibers  are rational curves.
Let $M\subset X(\r)$ be any connected component.
Then 
$$
M\sim N\ \#\ a\r\p^3 \ \#\ b(S^1\times S^2)\qtq{for some $a,b\geq 0$,}
$$
where  one of the following holds:
\begin{enumerate}
\item  $N$ is Seifert fibered over a  topological surface
(\ref{mult.fibr.def},
\ref{seifert.def}),
\item $N$ is the connected sum of  lens spaces (\ref{lens.defn}).
\end{enumerate}
\noindent  Assume in addition that $S$ is rational over $\c$. Then
one of three 
stronger restrictions holds:
\begin{enumerate}
\setcounter{enumi}{2}
\item  $N$ is Seifert fibered over a nonorientable surface
with at most 6 multiple fibers,
\item  $N$ is Seifert fibered over $S^2$ or $S^1\times S^1$ 
with at most 6 multiple fibers,
\item $N$ is the connected sum of at most 6 lens spaces.
\end{enumerate}
\end{thm}

\begin{rem}  I believe that (1--2) are sharp.

 The conclusions of (3--5) are, however, probably not
optimal. The computation of several examples suggests that the
following stronger assertions may hold:
\begin{enumerate}
\item In all  3  cases 6 can be replaced by 4. 
\item More precisely, if $m_i$
are the multiplicities of the multiple fibers
(or the lens spaces are $S^3/\z_{m_i}$) then
$\sum_i (1-\frac1{m_i})\leq 2$. 
\item If $N$ is Seifert fibered over   $S^1\times S^1$  then there are
no multiple fibers.
\end{enumerate}

Weaker results hold if $X(\r)$ is not orientable
(\ref{intro.nonor.case}).
\end{rem}

The map $f:X\to S$ may have finitely many 2-dimensional fibers.
The effect of these fibers can be investigated with the methods of
\cite{rat2}. In effect, \cite{rat2} shows that 
(\ref{orient.main.thm}) can be reduced to the case when we also assume
that every fiber of $f$ has dimension 1, at the price of allowing
certain isolated singularities on $X$.  Thus for the rest of this paper
I consider  morphisms $f:X\to S$ where every fiber of $f$ has dimension
1.
Two of the  resulting classes of morphisms is formalized in the next
definition:

\begin{defn} Let $X,Y$ be  normal varieties and $f:X\to Y$ a proper
morphism. $f$ is called a {\it rational curve fibration} if every fiber
of
$f$ has dimension 1 and the general fiber is a smooth rational curve.
(Special fibers may have several components. We do not   assume that
the irreducible components of the special fibers are rational, though
this is a consequence of the other assumptions \cite[I.3.17 and
II.2.2]{koll96}.)

$f$ is called a {\it conic bundle} if $Y$ is smooth and there is a
$\p^2$-bundle
$g:P\to Y$ and an embedding $j:X\into P$ such that
$f=g\circ j$ and every fiber of $f$ becomes a conic in the
corresponding $\p^2$. The conic  can be  smooth, a pair of lines
or a double line. (This definition is slightly more general
than the one given in \cite[Chap.I]{beau77} since we allow $X$ to be
singular.)

Note that $P=\proj_Yf_*\o_X(-K_{X/Y})$, thus $P$ is uniquely determined
by $X$.
\end{defn}

Assume now that $X,Y$ are real algebraic varieties and  $f$ is
defined over
$\r$ (see (\ref{what.is.rav?}) for the definitions). We can look at the
set of real points 
$X(\r)$ (resp. $Y(\r)$) 
and obtain a map of topological spaces
$f:X(\r)\to Y(\r)$. The preimage of a general point is the set of real
points of a smooth rational curve,  hence it is either empty or
homeomorphic to the circle $S^1$. Thus we conclude that $X(\r)$
contains a dense open set which is an $S^1$-bundle. Unfortunately,
every manifold has this property, and we need to understand the
behaviour of $f$ at every fiber in order to get useful topological
information about $X(\r)$. 

In Section 2 we prove the first general result in this direction:

\begin{thm}\label{intro.top.gen.cb.thm}
 Let $X$ be a projective 3--fold over $\r$ with isolated
singularities  and $f:X\to S$  a rational curve fibration. Assume
that
$X(\r)$ is an orientable PL  3--manifold.

Then there is a PL map $\tilde f:X(\r)\to F$ onto a surface with boundary
such that $\tilde f^{-1}(P)$ is a circle for every $P\in \inter F$ and
$\tilde f^{-1}(P)$ is a point for every $P\in \partial F$.
($\tilde f$ can be obtained as a small perturbation of $f|_{X(\r)}$
but it usually can not be chosen algebraic in a natural way.)
\end{thm}

\begin{rem} More generally, the above result also holds if we only
assume that
$\overline{X(\r)}$ is a PL 3--manifold.
($\overline{X(\r)}$ denotes the topological normalization of
$X(\r)$.  In our case, $X(\r)$ is a manifold except at finitely many
points where it is locally like the cone over a possibly disconnected
surface. The topological normalization   separates these local
components.) This is crucial for our applications.
\end{rem}

The topological  3--manifolds  
described in (\ref{intro.top.gen.cb.thm})
can be understood in terms of the
usual classification scheme of topological 3--manifolds. This is done
in section 3. 
 $\tilde f$  is an $S^1$-bundle over all but finitely many
points of $\inter F$ and it is exactly these finitely many points which
assert the greatest influence on the topology of $X(\r)$. 
In order to study these special fibers, we need a definition.

\begin{defn}\label{mult.fibr.def}
 Let $g:M\to F$ be a proper PL map of a 3--manifold to a surface
such that every fiber is a circle. Pick $P\in F$ and let
$P\in D_P\subset F$ be a small disc around $P$. Then
$g^{-1}(D_P)$ retracts to $g^{-1}(P)$. Let $P'\in D_P$ be a general
point. The retraction gives a map
$$
r(P',P):S^1\sim g^{-1}(P')\to g^{-1}(P)\sim S^1.
$$
The absolute value of the degree of $r(P',P)$ is independent of the
choices of $P'$, the retraction and the orientations.
It is  called the {\it multiplicity} of the fiber $g^{-1}(P)$ and it
is denoted by
$m_P(g)$. (A   detailed study of the local structure of such maps is
given in section 3.)

If $M$ is orientable, then $m_P(g)=1$ for all but finitely many points
$P\in F$. We say that $M$ is {\it Seifert fibered} over $F$ if
$m_P(g)\geq 1$ for every $P\in F$.

The multiple fibers with $m_P(g)>1$ are in close analogy
with multiple fibers of  elliptic surfaces.  
(Historically,  Seifert fibers
came first.) The $m_P(g)=0$ case 
does not have an  analog in complex algebraic geometry
since it leads to fibers which are homologous to zero.
\end{defn}

\begin{say}\label{int.say.conds}
 The proof of (\ref{intro.top.gen.cb.thm}) given in sections 
2 and 3 does not establish any link between the topological multiple
fibers and the local algebraic  nature of $f$. It turns out, however, 
that  a very close relationship can be established if we pose  
additional restrictions on $X$, for instance if we assume that $X$ is
smooth. For many applications this is, however, too restrictive.

One way to obtain rational curve fibrations in practice is through a
minimal model program.  \cite[1.9--11]{rat2}
suggests that we should therefore concentrate on 
rational curve fibration $f:X\to S$  which satisfy
 the following
properties. (These conditions may seem rather technical,
but they are very natural from the point of view of the 
minimal model program. See, for instance, \cite{koll87, KoMo98}.
Also, we only use these conditions to prove (\ref{intr.cart->cb.rem}),
after which one can forget about what terminal means.)
\begin{description}
\item[(\ref{int.say.conds}.1)]  $X$ is a real projective 3-fold   with
{\it terminal singularities} (\ref{term.sings.defn}) such that $K_X$ is
Cartier along
$X(\r)$ and
$\overline{X(\r)}$ is a PL 3--manifold.

\item[(\ref{int.say.conds}.2)] $-K_X$ is
$f$-ample.
\end{description}

For varieties over $\c$, the analogous  rational curve fibrations
have been investigated by many authors.  A systematic study of the
special fibers was begun in \cite{prok1, prok3, prok2}.  In these papers
a huge number of examples  is described and a partial classification is
given. A full description of the general case  seems still far away.
Usually, adding a real structure only complicates matters,
but the assumption that $K_X$ be Cartier along $X(\r)$
turns out to be extremely strong. 
(A similar situation happens for extremal contractions, as observed in
\cite{rat2}.) 

Section 4 uses the methods of \cite{Mori88} and Prokhorov
to obtain the following consequence:

\begin{prop}\label{intr.cart->cb.rem}
Assume that $f:X\to S$ satisfies the conditions 
(\ref{int.say.conds}.1--2).  Then there is a finite set $T\subset S$
such that
\begin{enumerate}
\item $f:X\setminus f^{-1}(T)\to S\setminus  T$ is a conic bundle,
\item $\red f^{-1}(s)$ is an irreducible, smooth and rational
curve for every real point $s\in T$. Moreover, $X$ has only isolated
hypersurface singularities along $\red f^{-1}(s)$ with the exception
of a pair of conjugate points $P,\bar P\in f^{-1}(s)$.
At these points $X$ is complex analytically isomorphic to the quotient
of a hypersurface singularity by a cyclic group $\z_m$. 
\end{enumerate}
\end{prop}

In section 5 we begin the study of the local structure of $f$ near
the special points $T$ of (\ref{intr.cart->cb.rem}).
First we establish that in a neighborhood of $s\in T$, $f$ can be
described as the quotient of a conic bundle by a cyclic group
(\ref{cover.thm}). This result is further developed to a complete
description in  (\ref{quotient.thm}).
This way we obtain a complete topological description
of $f:X(\r)\to S(\r)$ near the special points $T$.

In order to understand $f:X(\r)\to S(\r)$, we still need to analyze
the topology of conic bundles over $\r$.  A complete listing of the
cases seems rather difficult, but in section 7 we see that
one never obtains {\it multiple} Seifert fibers.
This gives the following precise relationship between the
algebraic and topological multiple fibers of $f$:
\end{say}

\begin{thm}\label{intro.local.seif.fib}
Let $X$ be a real algebraic threefold satisfying (\ref{int.say.conds}.1)
and $f:X\to S$ a rational curve fibration satisfying
(\ref{int.say.conds}.2).
Let $\tilde f:\overline{X(\r)}\to F$ be the PL map constructed in
(\ref{intro.top.gen.cb.thm}). Then  there is a one--to--one
correspondence between the sets:
\begin{enumerate}
\item  Multiple
Seifert fibers of $\tilde f$  of multiplicty $m(\tilde f)\geq 2$, and
\item  Multiple fibers of $f$ which are real analytically isomorphic to
the normal form
$$
\left(\p^2_{x:y:z}\times \a^2_{s,t}\right)/\z_m\supset
(x^2+y^2-z^2=0)/\z_m
\to \a^2_{s,t}/\z_m,
$$
where $\z_m$ acts by rotation with angle $2b\pi/m$ on $(s,t)$   for some
$(b,m)=1$, it fixes
$z$ and acts by rotation with angle $2\pi/m$ on $(x,y)$.
\end{enumerate}
\end{thm}

The most interesting piece of information truns out to be the
description of the singularity occurring on $S$.
The quotient $\a^2_{s,t}/\z_m$ is  real analytically isomorphic to
the singularity  denoted by $A^+_{m+1}$
given by equation $(u^2+v^2-w^m=0)$ (cf.\ (\ref{dv.defn})).  It turns out
that it is quite hard for a real surface to contain such a singularity.
In section 9 we prove that if
$S$ is a real surface which is rational over $\c$
then
a connected component  of $\overline{S(\r)}$
   contains at most 6 singularities of type
$A^+$. (By contrast, a rational surface  can contain arbitrary many
singularities of type  $(u^2-v^2-w^m=0)$.)

All these results are assembled in section 8 to
 obtain the proof of the main theorems.

Finally, section 10 contains some examples of 3--manifolds
which can be realized by rational curve fibrations over $\r$.

\begin{say}[The nonorientable case]\label{intro.nonor.case}{\ }

The situation becomes more complicated if we do not assume that
$X(\r)$ is orientable. Because of \cite{rat2}, in this case it is more
natural to consider the case when $X$ satisfies
(\ref{int.say.conds}.1) and $f:X\to S$ is a rational curve fibration.

If $S$ is a (not necessarily rational)   real algebraic
surface then $\overline{X(\r)}$ is a connected sum of lens spaces
with a 3--manifold where all the pieces of the 
Jaco--Johannson--Shalen decomposition (cf.\ \cite[p.483]{Scott83})
are Seifert fibered. There are further restrictions but they are
somewhat complicated. Also, I do not know how sharp these results are.
See (\ref{nonorient.top.thm}) for details.

If $S_{\c}$ is rational then we should get very few cases, but  again
I do not have a reasonably sharp answer.
\end{say}

\begin{defn}\label{what.is.rav?}
 By a {\it real algebraic variety} I mean a 
variety given by real equations, as defined in most algebraic
geometry books (see, for instance, 
\cite{Shafarevich72, Hartshorne77}). 
This is consistent with the usage of \cite{Silhol89} but
 is  different from the
definition  of  \cite{BCR87}  which
essentially considers only the germ of $X$ along its real points. 
In many cases  the two variants can be used interchangeably, but
in this paper it is crucial to use the first one.

If $X$ is a real algebraic
variety then      $X(\r)$ denotes the set of real
points of $X$   as a topological space
and $X(\c)$ denotes the set of complex points as a complex space. 
$X_{\c}$ denotes the corresponding complex variety
(same equations as for $X$ but we pretend to be over $\c$).

For all practical purposes we can identify $X$ with the pair
($X(\c)$, complex conjugation) (cf.\ \cite[Sec.I.1]{Silhol89}).

A property of $X$ always refers to the variety $X$. Thus, for
instance, 
$X$ is   smooth iff it is   smooth at all complex points, not just at
its real points.  I use the adjective {\it ``geometrically"} to denote
properties of the complex variety $X_{\c}$. 
\end{defn}

\begin{say}[Piecewise linear 3--manifolds]\label{1.basic.top.facts}{\ }

In this paper I usually work with piecewise linear  manifolds
(see \cite{RoSa82} for an introduction).
Everey real algebraic variety carries a natural PL structure
(cf.\ \cite[Sec.9.2]{BCR87}).

 In dimension 3 every
compact topological 3--manifold carries a unique PL--manifold
structure (cf.\ \cite[Sec.\ 36]{Moise77}) and  a PL--structure behaves
very much like a differentiable  structure. For instance, let $M^3$ be a
PL 3--manifold,
$N$ a compact PL--manifold of  dimension 1 or 2 and $g:N\into M$ a
PL--embedding. Then a suitable open  neighborhood of $g(N)$ is
PL--homeomorphic to a real vector bundle over $N$  (cf.\ \cite[Secs.\
24 and 26]{Moise77}).  (The technical definition of these is given by
the notion of  {\it regular neighborhoods}, see \cite[Chap.3]{RoSa82}.) 
If $f:M\to N$ is a PL--map and $X\subset N$ a compact subcomplex then
there is a regular neigborhood $X\subset U\subset N$ such that
$f^{-1}(U)$ 
is a regular neigborhood of $f^{-1}(X)\subset M$ 
 (cf.\ \cite[2.14]{RoSa82}). 
\end{say}

\begin{ack}  I   thank S. Gersten, M. Kapovich and B. Kleinert for
answering my numerous questions about  3-manifold topology.  
Partial financial support was provided by  the NSF under grant number 
DMS-9622394. 
\end{ack}

\section{The topology of rational curve fibrations over $\r$}

The aim of this section is to provide an argument
describing the topology of  rational curve fibrations  over $\r$ under
rather general assumptions. A weakness of the result is that it does not
provide any connection between the algebraic structure of the singular
fibers and the topology over $\r$.  The precise relationship is worked
out in subsequent sections.

\begin{thm}\label{top.gen.cb.thm}
 Let $X$ be a projective 3--fold over $\r$ with isolated
singularities only and $f:X\to S$ a rational curve fibration. Assume that
$\overline{X(\r)}$ is an orientable  3--manifold.

Let  $M\subset \overline{X(\r)}$   be any connected component.
Then there are disjoint solid tori
$S^1\times D^2\sim T_i\subset M $ such that
$M\setminus\cup_i\inter T_i$ is an $S^1$-bundle over a  compact
surface with boundary.
\end{thm}

We need the description of rational curve fibrations over curves.
These are discussed in 
 \cite[Secs.\ II.6, VI.3]{Silhol89} and \cite[1.8]{ras}.

\begin{lem}\label{rcfibr.surf.lem}
 Let $Y$ be a smooth projective surface over $\r$  and
$f:Y\to A$ a rational curve fibration.  Every point $0\in A(\r)$ has a
neighborhood  $0\in U\subset  A(\c)$ such that $f^{-1}(U)\to U$  is real
analytically equivalent to one of the following normal forms:
(In all 4 cases $f$ is the
second projection and $D$ is the unit disc in $\a^1_s$.)
\begin{enumerate}
\item ($\p^1$-bundle) 
$(x^2+y^2-z^2)\subset \p^2_{(x:y:z)}\times D$,
\item (empty fibers) 
$(x^2+y^2+z^2)\subset \p^2_{(x:y:z)}\times D$,
\item (collapsed end)
$(x^2+y^2+sz^2)\subset \p^2_{(x:y:z)}\times D$,
\item (blown up cases) obtained from one of the above by repeatedly
blowing up real points and conjugate pairs  of complex points.\qed
\end{enumerate}
\end{lem}

\begin{say}[Proof of (\ref{top.gen.cb.thm})]{\ }

We may assume that $X$ and $S$ are normal.

Let $A\subset S$ be a general hyperplane
section and set $Y=f^{-1}(A)$.  $Y(\r)$ is an orientable 2--manifold.
 Since
$Y(\r)$ is orientable, we can only blow up conjugate pairs of complex
points in (\ref{rcfibr.surf.lem}.4)   (cf.\ \cite[Sec.II.6]{Silhol89},
\cite[2.2]{ras}). Thus
$Y(\r)\to A(\r)$  is locally (on $A(\r)$) homeomorphic to one of three
normal forms:
\begin{enumerate}
\item $S^1$-bundle
\item empty fibers
\item $(x^2+y^2+sz^2)\subset \r\p^2_{(x:y:z)}\times D$,
with projecion to $\a^1_s$. 
\end{enumerate}
 Let $B\subset S$ be
the  locus of singular fibers and 
set $U=f(X(\r))$. $U\subset S(\r)$ is a semi algebraic subset.
Aside from a finite set $T\subset S(\r)$,
we have the following description:
\begin{enumerate}
\item $\partial U=B(\r)$,
\item $f:X(\r)\to U$ is an $S^1$-bundle over $\inter U$,
\item each fiber over $\partial U$ is a single point.
\end{enumerate}
We still have to determine the local structure of $X(\r)\to U$ near the
points in $T$. Let $P\in T$ be a point. 
$f^{-1}(P)(\r)$ is a real algebraic curve, hence a union of some copies
of $S^1$. The normalization  $\overline{X(\r)}\to X(\r)$ can break
apart some of the circles, thus the preimage of $f^{-1}(P)(\r)$ in 
$\overline{X(\r)}$ is a
1--complex; let $D_i(P)$ be its connected components.

The boundary of a regular neighborhood of
$D_i(P)$ in $\overline{X(\r)}$ is   connected, hence the map
$\overline{X(\r)}\to U$ factors through the normalization $\bar U\to U$.
We denote it by  $\bar f:\overline{X(\r)} \to \bar U$.
We are in the following situation:
\begin{enumerate}
\item $\overline{X(\r)}$ is a PL 3-manifold and 
$\bar U$ is a PL 2--manifold with boundary.
\item There is a finite set $T\subset U$ such that
\begin{enumerate}
\item $f^{-1}(P)$ is  1--complex for $P\in T$,
\item $f$ is an $S^1$-bundle over $\inter U\setminus T$, and
\item $f^{-1}(\partial U\setminus T)\to (\partial U\setminus T)$ is a PL
homeomorphism.
\end{enumerate}
\end{enumerate}

Pick $P\in T$. 
If $h^1(D_i(P),\q)=g_i$,  then the boundary of a regular
neighborhood of
$D_i(P)$ in $\overline{X(\r)}$ is an orientable surface of genus
$g_i$ (cf.\ \cite[2.4]{Hempel76}). 
On the other hand, this boundary is an $S^1$-bundle over $S^1$
if $P\in \inter U$   and $S^2$ if
$P\in \partial U$ by (\ref{1.basic.top.facts}). Thus $\bar f^{-1}(P)$ is
collapsible to an
$S^1$ in the first case and to a point in the second case. We can do
these collapsings in 
$\overline{X(\r)}$.  This creates a 3--manifold $\overline{X(\r)}^*$
which is PL homeomorphic to $\overline{X(\r)}$ (cf.\
\cite[3.27]{RoSa82}).  We replace $\overline{X(\r)}$ with a connected
component $M$ of $\overline{X(\r)}^*$ to obtain 
$g:M\to V$, satisfying the following conditions:
\begin{enumerate}
\item $M$ is a compact PL 3-manifold and 
$V$  a PL 2--manifold with boundary.
\item There is a finite set $T\subset \inter V$ such that
\begin{enumerate}
\item $f^{-1}(P)\sim S^1$   for $P\in T$, 
\item $f$ is an $S^1$-bundle over $\inter U\setminus T$, and
\item $f^{-1}(\partial V)\to \partial V$ is a PL homeomorphism.
\end{enumerate}
\end{enumerate} 

Let $F^0\subset V$ be obtained by removing small open discs $U_P$
around each $P\in T$ and a small annulus  $U_C$ along each boundary
component $C\subset \partial V$. Then $M^0:=f^{-1}(F^0)$ is an
$S^1$-bundle over
$F^0$. Moreover, each $f^{-1}(\bar U_P)$ is  a
regular neighborhood of $f^{-1}(P)\sim S^1$
and each $f^{-1}(\bar U_C)$ is  a
regular neighborhood of $f^{-1}(C)\sim S^1$. By the orientability
of $M$, the $f^{-1}(\bar U_P)$ and the $f^{-1}(\bar U_C)$ are 
 solid tori.
\qed
\end{say}

\section{Surgery on circle bundles}

The aim of this section is to give a more detailed description
of the 3--manifolds that appear in the study of real conic bundles.
For the purposes of (\ref{top.gen.cb.thm}) we need only the orientable
case, but later on we encounter 
nonorientable examples as well.

\begin{defn}[Lens spaces]\label{lens.defn}{\ }
For relatively prime $0<q<p$ consider the
action
$(x,y)\mapsto (e^{2\pi i/p}x, e^{2\pi iq/p}y)$ on the unit sphere
$S^3\sim (|x^2|+|y^2|=1)\subset \c^2$. The quotient is a 3--manifold,
called the {\it lens space}  $L_{p,q}$. 

Another way to obtain lens spaces is to glue two solid tori together.
The result is   a lens space,  $S^3$ or $S^1\times S^2$. 
Sometimes one writes $L_{1,0}=S^3$ and $L_{0,1}=S^1\times S^2$.
See, for instance,
\cite[p.20]{Hempel76}.
\end{defn}

\begin{say}\label{def.tor-stor}
[Torus and solid tous]{\ }

As a general reference, see \cite[Chap.2]{Rolfsen76}.

Let $S^1$ be the unit circle $(|u|=1)\subset \c_u$ and $D^2$ the unit
disc
$(|z|\leq 1)\subset \c_z$. 
$S^1\times D^2$ is called the {\it solid torus}. Its boundary
$$
\partial (S^1\times D^2)=(|u|=1)\times (|z|=1)\sim S^1\times S^1
$$
is a {\it torus}. Up to isotopy, the homeomorphisms of a torus are given
by
$$
(u,z)\mapsto (u^az^b,u^cz^d),\qtq{where}
\left(
\begin{array}{cc}
a& b\\
c&d
\end{array}
\right)
\in GL(2,\z).
$$
Up to isotopy, the homeomorphisms of a solid torus are given by
$$
(u,z)\mapsto (u^{\pm 1},u^mz)\qtq{or}
(u,z)\mapsto (u^{\pm 1},u^m\bar z),\qtq{where} m\in\z.
$$
Let $g:S^1\times S^1\to S^1$ be an $S^1$-bundle. Up to isotopy,
it can be written as
$$
g_{c,d}:(u,z)\mapsto u^cz^d, \qtq{where} (c,d)=1.
$$
Given a solid torus $S^1\times D^2$, let 
$g_{c,d}:\partial (S^1\times D^2)\to S^1$ be an $S^1$-bundle.
This can be extended to a  map
$$
G_{c,d}:S^1\times D^2\to D^2 \qtq{by}
(u,z)\mapsto
\left\{
\begin{array}{lll}
u^cz^d&\mbox{if}& d> 0,\\
u^c\bar z^{(-d)}&\mbox{if}& d< 0,\\
u^c|z| &\mbox{if}& d= 0.
\end{array}
\right.
$$
The fiber $G_{c,d}^{-1}(1)$ is parametrized by
$$
\phi:S^1\to S^1\times D^2\qtq{given as} t\mapsto (t^d,t^{-c}).
$$
($\phi$ is  injective since $(d,c)=1$.) Composing $\phi$ with the
first projection $S^1\times D^2\to S^1$ we obtain a map of degree $m$
from $S^1$ to $S^1$. Thus the fiber of $G_{c,d}$ over the origin
has multiplicity $|d|$ (\ref{mult.fibr.def}).

Applying a homeomorphism of the solid torus, we can get
every $G_{c,d}$ to
one of the following normal forms:
\begin{enumerate}
\item 
For a pair of integers $c,d$ satisfying  $0\leq c<d$ and
$(c,d)=1$, define
$$
f_{c,d}:S^1\times D^2\to D^2\qtq{by}
 f_{c,d}(u,z)=u^cz^d.
$$
$f_{c,d}$ restricts to a fiber bundle
$S^1\times (D^2\setminus \{0\})\to D^2\setminus \{0\}$.
The fiber of $f_{c,d}$ over the origin is still $S^1$, but 
$f_{c,d}^{-1}(0)$ has
multiplicity $d$. 
\item $f_{1,0}(u,z)=u\cdot |z|$.  In this case the fibers
$f_{1,0}^{-1}(t)$ for $t\neq 0$ are null homotopic in 
the solid torus. 

Instead of $f_{1,0}$, I prefer to use the map
$$
f_0:(u,z)\mapsto u\frac{1+|z|}{2}\in \{t\vert 1/2\leq |t|\leq 1\}\subset
\c_t.
$$
The image is a closed annulus  $A_{1/2,1}$ with boundary circles
$C_{1/2}$ and $C_1$. 
$f_0$ is an $S^1$-bundle over the annulus $A_{1/2,1}\setminus C_{1/2}$
and $f_0^{-1}(C_{1/2})\to C_{1/2}$ is a homeomorphism.

$f_{1,0}$ is obtained by composing $f_0$ with $t\mapsto
2t(1-\frac{1}{2|t|})$ which contracts the circle $C_{1/2}$ to a point. 
\end{enumerate}
\end{say}

\begin{say}\label{def.kb-skb}
[Klein bottle and solid Klein bottle]{\ }

On the solid torus above, consider the map $(u,z)\mapsto (u^{-1},\bar
z)$. This is a fixed point free and orientation reversing ivolution.
The quotient space is called a {\it solid Klein bottle}. Its boundary is
a {\it Klein bottle}.

Among the homeomorphisms of the torus, only the maps
 $(u,z)\mapsto (u^{\pm 1}, z^{\pm 1})$ descend to the Klein bottle,
and these are all the homeomorphisms of the Klein bottle modulo isotopy.
All of these extend to the solid Klein bottle as 
$$
(u,z)\mapsto (u^{\pm 1}, z)\qtq{or} (u,z)\mapsto (u^{\pm 1}, \bar z).
$$
Thus there is a unique way to glue a solid Klein bottle to
 a manifold along a  boundary component which is a  Klein bottle.

Given an $S^1$-bundle structure on the boundary Klein bottle, the
second normal form in (\ref{def.tor-stor}) gives the extension
$$
f_0:(u,z)\mapsto u^2\frac{1+|z|}{2}\in \{t\vert 1/2\leq |t|\leq
1\}\subset
\c_t.
$$
\end{say}

\begin{defn}[Seifert fiber spaces]\label{seifert.def}{\ }

 A {\it Seifert fibered} 3--manifold is a
proper morphism of a 3--manifold to a surface $f:M\to S$
such that every point $s\in S$ has a neighborhood $s\in D_s\subset S$
such that the pair $f^{-1}(\bar D_s)\to \bar D_s$ is fiber preserving
homeomorphic to one of the  normal forms $f_{c,d}$  defined 
in (\ref{def.tor-stor}.1)
for some $c,d$ satisfying  $0\leq c<d$ and
$(c,d)=1$.

There are other ways of describing the local models.
Consider the map
$$
h_{c,d}:[0,1]\times D^2 \to S^1\times D^2
\qtq{given by} (t,z)\mapsto (e^{2\pi it}, ze^{-2\pi ict/d}).
$$
Under this map we can view $S^1\times D^2$ as 
$[0,1]\times D^2$ with the two ends $\{0\}\times D^2$ and 
$\{1\}\times D^2$ identified by clockwise rotation by $2\pi c/d$.
Moreover, each $[0,1]\times \{z\}$ maps to a single fiber of
$f_{c,d}$. This shows that the above definition is equivalent with the
usual one (cf.\ \cite{Scott83}).

Finally, we can 
consider the action of $\z_d$ on $S^1\times D^2$ which is
rotation by $2\pi c/d$ on $S^1$ and by $2\pi/d$ on $D^2$. The quotient
$(S^1\times D^2)/\z_d\to D^2/\z_d\sim D^2$ is Seifert fibered
with multiplicity $d$.

(In several papers, one   nonorientable  local model is also
allowed. For
$p=2$ we can act on $D^2$ by reflection. Then $D^2/\z_2$ is a surface
with boundary and $f$ is not a fiber bundle near the boundary. This case
does not come up for us. An algebraic model of this situation is
given by
$$
X:=(z^2+s(x^2+y^2)=0)\subset \p^2_{x,y,z}\times \a^2_{s,t},
$$
where $f$ is the second projection. Notice that
$X$ is singular along the curves $(z=s=x\pm \sqrt{-1}y=0)$. This is a
good example to show how complex singularities influence the behaviour
of the real points.)
\end{defn}

\begin{thm}\label{top.sf-ls.thm}
 Let $F^0$ be a compact PL surface with boundary
and $F\supset F^0$ the surface obtained from $F^0$ by attaching discs
to all boundary components. Let $M^0\to F^0$ be an $S^1$-bundle (so
its boundary components are tori and Klein bottles). Let $M$
be any 3--manifold obtained from
$M^0$ by  attaching solid tori or Klein bottles
to all boundary components.
Then one of the following holds:
\begin{enumerate}
\item  $M$ is Seifert fibered over $F$, or
\item  $M$ is the connected sum of lens spaces,  $S^1\times S^2$  and
 $S^1\tilde{\times} S^2$. 
\end{enumerate}
\end{thm}

\begin{rem} In the orientable case, (\ref{top.sf-ls.thm})
is equivalent to the following statement.

If $M$ is obtained from an $S^1$-bundle by surgery along fibers then
either $M$ is Seifert fibered or $M$ is the connected sum of lens
spaces and  $S^1\times S^2$. 
\end{rem}

\begin{complement}\label{top.sf-ls.compl}
 Notation as in (\ref{top.sf-ls.thm}). The gluing data
determine the structure of $M$ as follows.

Let $A_1,\dots, A_r$ be the boundary components of $M^0$ homeomorphic to 
a torus and $T_1,\dots, T_r$ the corresponding solid tori attached by
$\phi_i:\partial T_i\sim A_i$.
On each $T_i$  choose   coordiantes $(u_i,z_i)$ as in
(\ref{def.tor-stor}). The $S^1$-bundle structure of $M^0$ induces
$S^1$-bundles $p_i:A_i\to S^1$ for every $i$. We may assume that
$p_i\circ \phi_i=u_i^{c_i}z_i^{d_i}$,
where $c_i,d_i$ satisfy  the usual conditions $0\leq c_i<d_i$ and
$(c_i,d_i)=1$ or $c_i=1, d_i=0$.

Let $B_1,\dots, B_s$ be the Klein bottle boundary components.
\begin{enumerate}
\item If $s=0$ and $d_i>0$ for every $i$ then $M$ is Seifert fibered
over $F$ with mulitple fibers of multiplicity $d_1,\dots, d_r$.
\item If $s>0$ or $d_i=0$ for some $i$, then 

\noindent $M\sim L_{d_1,c_1}\ \#\ \cdots\ \#L_{d_r,c_r}\ \#\ a(S^1\times
S^2)
\ \#\ b(S^1\tilde{\times} S^2)$

\noindent for some $a,b\geq 0$. 
\end{enumerate}
\end{complement}

Proof. The local models of attaching a solid torus or solid Klein bottle
are described in (\ref{def.tor-stor}) and (\ref{def.kb-skb}).
We obtain that there are
\begin{enumerate}
\item  a surface $F^0\subset F^1\subset F$, obtained from $F^0$ by
attaching discs or annuli, and
\item  a proper map $g:M\to F^1$ such that
$g^{-1}(\inter F^1)\to \inter F^1$ is Seifert fibered and
$g^{-1}(\partial F^1)\to \partial F^1$ is a homeomorphism.
\end{enumerate}

If $\partial F^1=\emptyset$ then we are done.

Assume next that $\partial F^1\neq \emptyset$ and let $L\subset F^1$ be a
simple path starting and ending in $\partial V$. Then $g^{-1}(L)\sim
S^2$, hence cutting along $g^{-1}(L)$ corresponds to connected sum
decomposition (if $L$ separates $F^1$) or to removing a 1--handle
(if $L$ does not separate $F^1$), cf.\ \cite[3.8]{Hempel76}. After
repeated cuts we are reduced to the situation when 
each of the pieces $\cup V_i=F^1$
 is a disc 
containing at most 1 multiple fiber
and
$M$ is the connected sum of the corresponding manifolds $M_i$
(possibly with further 1--handles attached).
$g_i:M_i\to V_i$ 
has a unique multiple fiber in   $\inter V_i$ at  a point $P_i$
and $g^{-1}(\partial V_i)\sim \partial V_i\sim S^1$. A
regular neighborhood of
$g^{-1}(P_i)$ is a solid torus and so is a regular neighborhood of
$g^{-1}(\partial V_i)$. Thus $M_i$ is obtained  by gluing two
solid tori together, hence $M_i$ is   a lens space, $S^3$ or   $S^1\times
S^2$, depending on  $d_i,c_i$ (\ref{lens.defn}). \qed

\section{General results on rational curve fibrations over $\c$}

 Let $X$ be a 3--fold with terminal singularities and
$f:X\to S$ a rational curve fibration such that $-K_X$ is $f$-ample. A
considerable effort has been spent on trying to understand the local
structure of $f$ (see
\cite{prok1, prok3, prok2}), but so far we do not have a complete
description. The examples of Prokhorov suggest that there are many cases.
The aim of this section is to prove   a technical
result on rational curve fibrations (\ref{cb.gen.prop}). The  proof is 
quite easy but it uses some of the
  machinery of extremal neighborhoods developed in
\cite{Mori88}. This result is used only through (\ref{cart->cb.rem}),
so readers more interested in the topological aspects can skip the rest
of the secion.

\begin{say}[Terminal singularities]\label{term.sings.defn}{\ }

For a precise definition of terminal singularities the reader is
referred to \cite{Reid85, KoMo98}. We need only the following
consequence of their classification:

Every  3--dimensional terminal
singularity is either
\begin{enumerate}
\item an isolated hypersurface singularity $(0\in X)\subset \c^4$,
(these are the cases when $K_X$ is Cartier), or
\item  a quotient of an  isolated hypersurface singularity by a
$\z_m$-action with an isolated fixed point. 
\end{enumerate}
\noindent The value of $m$ is
called the {\it index} of the singularity.

The classification of 3--dimensional terminal
singularities over $\r$ is given in \cite{rat1}. We do not need it but
in some cases I use their local topological description 
developed in \cite[Secs.\  4--5]{rat1}. 
\end{say}

\begin{prop}\label{cb.gen.prop}
 Let $X$ be a 3--fold with terminal singularities over $\c$ and
$f:X\to S$   a rational curve fibration 
such that $-K_X$ is $f$-ample. Let 
$0\in S$ be a closed point. 
Then  one of the following   alternatives holds:
\begin{enumerate}
\item $K_X$ is Cartier along $f^{-1}(0)$, or
\item 
$K_X$ is not Cartier
at all singular points of $\red f^{-1}(0)$.
\end{enumerate}
\end{prop}

We are mainly interested in the following consequence
for  rational curve fibrations over $\r$:

\begin{cor}\label{cart->cb.rem}
Let  $X$ be a real projective 3-fold with terminal
singularities such that $K_X$ is Cartier along $X(\r)$.
Let    $f:X\to S$ be a    rational curve fibration over $\r$ such that
$-K_X$ is
$f$-ample. For every point  $s\in S(\r)$, one of the following holds: 
\begin{enumerate}
\item $S$ is smooth at $s$ and $f$
is  a conic bundle over a neighborhood of $s$, or
\item $\red f^{-1}(s)$ is an irreducible, smooth and rational
curve and $K_X$ is not Cartier precisely at a conjugate pair of 
complex points of $\red f^{-1}(s)$.
\end{enumerate}
\end{cor}

Proof.  If $K_X$ is Cartier along $f^{-1}(s)$ then $S$ is smooth at 
$0$ and $f$ is  a conic bundle over a neighborhood of $s$ by
\cite[Thm7]{Cutkosky88}. This gives  (\ref{cart->cb.rem}.1). 

If $\red f^{-1}(s)$ is an irreducible, smooth and rational
curve  then 
$X$ can have at most 3 singular points along $\red f^{-1}(s)$
by \cite[0.4.13.1]{Mori88} and \cite[2.1]{prok1}.  Thus either
$K_X$ is Cartier along $\red f^{-1}(s)$ or there is a unique
conjugate pair of 
complex points of $\red f^{-1}(s)$ where $K_X$ is not Cartier.
Thus we have (\ref{cart->cb.rem}.2).

In all other cases, $\red f^{-1}(s)_{\c}$ is a  tree of smooth rational
curves  by (\ref{conic.fib.defs}) and $K_X$ is not Cartier
at all singular points of $\red f^{-1}(s)_{\c}$ by 
(\ref{cb.gen.prop}).   $\red f^{-1}(s)_{\c}$ can not have  a real
singular point by assumption, hence it contains a real irreducible
component $D\subset \red f^{-1}(s)_{\c}$. As in 
(\ref{conic.fib.defs}),  $D$ can be contracted in an extremal
contraction. By \cite[7.2]{rat2} this is only possible if   $K_X$ is 
 Cartier along
$D$, a contradiction.
\qed

\begin{say}\label{conic.fib.defs}
[Rational curve fibrations and extremal neighborhoods]{\ }

Let $X$ be a 3--fold with terminal singularities and $g:X\to Y$
 a birational morphism. Assume that $-K_X$ is $g$-ample
and $g^{-1}(y)$ is a curve for some $y\in Y$. The germ of $X$ along 
$g^{-1}(y)$ is called an {\it extremal curve neighborhood}
(cf.\ \cite[p.549]{KoMo92}). Extremal curve neighborhoods have been
studied in detail in \cite{Mori88, KoMo92}. A key point
in their study are the vanishings $R^1g_*\o_X=0$ and
$R^1g_*\omega_X=0$. The first of these still holds for rational curve
fibrations by the general Kodaira vanishing theorem
(cf.\ \cite[Sec.2.5]{KoMo98}) but not the second. 

By an easy argument (cf.\ \cite[1.2.1]{Mori88}), this implies the
following: If $\o_X\onto Q$ is  any quotient supported on $f^{-1}(0)$
then $H^1(X,Q)=0$. In particular, $H^1(\o_{\red f^{-1}(0)})=0$.
Any such curve is a tree of smooth rational curves (cf.\
\cite[1.3]{Mori88}).

For many results in \cite{Mori88}, these are the key vanishings
 needed. In practice this means that
many of the results for extremal curve neighborhoods still hold for
rational curve fibrations, maybe in slightly modified form.
This approach has been carried out   by \cite{prok3}.

In some cases one can directly apply the results on 
extremal curve neighborhoods to rational curve fibrations through the
following trick:

 Pick $0\in S$ and let $\cup_iC_i=C=\red f^{-1}(0)$ be the
irreducible components. If we replace $S$
by a small analytic neighborhood of $0$ then there are divisors
$H_i\subset X$ such that $H_i$ intersects $C$ only at a point of $C_i$.
Let $D\subsetneq C$ be any closed subcurve and $H_D$ the sum of those
$H_i$ which are disjoint from $D$.
 By the relative base point free
theorem (cf.\ \cite[Sec.3.6]{KoMo98}), a mutiple of $H$ gives a proper
birational morphism
$p_D:X\to Y$. 
$p_D$ contracts $D$ to a point ($p_D$ may also contract other curves
in other fibers of $f$).
$-K_X$ is $p_D$-ample, thus $D$ determines an extremal curve
neighborhood in
$X$. In particular, all the results
of \cite{KoMo92} describing extremal neighborhoods apply to $D$.
This of course does not work if $D=C$.
\end{say}

\begin{say}\label{adj.formula}[Adjunction formula and $i_P(1)$]{\ }

If $Y$ is a smooth 3--fold and $C\subset Y$ a smooth curve 
with ideal sheaf $I_C$, then the
adjunction formula says that 
$$
2g(C)-2=(C\cdot K_Y)-\deg_C(I_C/I_C^2).
$$
If $Y$ is singular but $C$ smooth,   a similar formula was
developed in \cite{Mori88}, which also involves a correction term
depending on the singularities of $X$ along $C$. The two corrections
coming from $(C\cdot K_Y)$ and from $\deg_C(I_C/I_C^2)$
are mixed together  into a single number $i_{P,C}(1)$
(denoted by $i_P(1)$ in \cite{Mori88} and by $i_P$ in \cite{prok1}). 
The definition is not crucial for us, and we need only two properties:
\begin{enumerate}
\item $i_{P,C}(1)$ is a nonnegative integer which     
is zero iff
$P$ is a smooth point of $Y$ \cite[2.15]{Mori88},
\item If $C\cong \p^1$ gives  an extremal curve neighborhood then
$\deg_{C}I_C/I_C^{(2)}=1-\sum_{P\in C}i_{P,C}(1)$ \cite[2.3.2]{Mori88}.
\end{enumerate}
\end{say}

\begin{say}[Proof of (\ref{cb.gen.prop})]{\ }

  Write $f^{-1}(0)=\cup_iC_i$.
$\cup_iC_i$ is a tree of smooth rational curves, hence all
the singular points of 
$\red f^{-1}(0)$ are at intersection points of the irreducible
components. Let $P\in C_i\cap C_j$ be such a point and set
$D=C_i\cup C_j$. 

If $D\neq \red f^{-1}(0)$ then $D$ determines an extremal neighborhood,
so  $K_X$ is not Cartier at 
$P$ by
\cite[1.15]{Mori88}. Thus (\ref{cb.gen.prop}.2) holds if $f^{-1}(0)$ has
at least 3 irreducible components. (\ref{cb.gen.prop}.2)   holds
vacuously if 
 $f^{-1}(0)$ has only one irreducible component.
We are left with the case when $f^{-1}(0)=C_1\cup C_2$ has two
irreducible components. They have a unique intersection point $Q$.
We are done if $K_X$ is not Cartier at $Q$.

Assume next that $K_X$ is  Cartier  along $C_1$. Let $p:X\to Y$ be the
contraction of $C_1$. By \cite[Thm.4]{Cutkosky88},
$p$ is a divisorial contraction, $Y$ again
has terminal singularities and $K_X=p^*K_Y+E$
where $E$ is a $p$-exceptional Cartier divisor.
 Let  $g:Y\to S$ be the induced rational curve fibration.
A general fiber of
$g$ is algebraically equivalent to a  multiple $n\cdot p(C_2)$
of $p(C_2)$. If $n=1$ then $g$ is a $\p^1$-bundle by
\cite[II.2.28]{koll96}, thus  $K_X$ is Cartier along $C_2$ as well and we
are in case (\ref{cb.gen.prop}.1).
Otherwise  $(K_Y\cdot 
p(C_2))\geq -2/n\geq -1$ and so
$$
(K_X\cdot C_2)\geq (K_Y\cdot  p(C_2))+(E\cdot C_2)\geq -1+1=0,
$$
 a contradiction. 

We are left with the case when $K_X$ is  Cartier at $Q$ but
not Cartier  along some point on either $C_i$.
Let $I_i\subset \o_X$ be the ideal sheaf of
$C_i$. Consider the exact sequence
$$
\begin{array}{ccccc}
0&\to & H^0(\o_X/(I_1^{(2)}\cap I_2^{(2)}))
&\to & H^0(\o_X/I_1^{(2)}) + H^0(\o_X/I_2^{(2)})\\
&\to & H^0(\o_X/(I_1^{(2)}+ I_2^{(2)}))
&\to & H^1(\o_X/(I_1^{(2)}\cap I_2^{(2)})).
\end{array}
$$
$H^1(\o_X/(I_1^{(2)}\cap I_2^{(2)}))=0=H^1(\o_X/I_i^{(2)})$ 
as we noted in (\ref{conic.fib.defs}). By (\ref{adj.formula}.2)
and Riemann--Roch on $C_i$,
$$
h^0(\o_X/I_i^{(2)})=\chi(\o_X/I_i^{(2)})=3+\deg_{C_i}(I_i/I_i^{(2)})
=4-\sum_{P\in C_i}i_{P,C_i}(1).
$$
$\o_X/(I_1^{(2)}+ I_2^{(2)})$ is supported at $Q$.
In order to compute it, we consider
 two cases:
\begin{enumerate}
\item ($X$ is smooth at $Q$) In suitable local coordinates
$I_1=(x,y)$ and $I_2=(y,z)$. Thus
$\o_X/(I_1^{(2)}+ I_2^{(2)})$ is spanned by $1,x,y,z,xz$ and
 $h^0(\o_X/(I_1^{(2)}+ I_2^{(2)}))=5$.
\item ($X$ is singular at $Q$) In suitable local coordinates
$I_1=(x,y,z)$ and $I_2=(y,z,t)$. Thus
$\o_X/(I_1^{(2)}+ I_2^{(2)})$ is spanned by $1,x,y,z,t,xt$ and if the
equation of $X$ involves $xt$, then $xt$ is zero in the quotient.
Thus $h^0(\o_X/(I_1^{(2)}+ I_2^{(2)}))\geq 5$.
\end{enumerate}
Comparing these with the above exact sequence we conclude that
$$
8-\sum_{P\in C_1}i_{P,C_i}(1)-\sum_{P\in C_2}i_{P,C_i}(1)\geq 1+5.
$$
If $X$ is not smooth at $P\in C_i$ then $i_{P,C_i}(1)\geq 1$ by
(\ref{adj.formula}.1).
Since $K_X$ is not Cartier  along either $C_i$, we conclude that 
there is a
unique  point $P_i\in C_i$ such that  $i_{P_i,C_i}(1)= 1$. 

Each $C_i$ defines an extremal curve neighborhood in $X$,
and we are in cases (2.2.1) or (2.2.1') in the list
\cite[2.2]{KoMo92}.  (In the cases (2.2.2) or (2.2.2') the
singular point has
$i_{P,C}(1)=2$ by \cite[6.5]{Mori88} and \cite[p.549]{KoMo92}.
In the cases (2.2.3), (2.2.3') and (2.2.4) there are at least 2
singular points.)
From this we conclude that  $-K_X$ has a good member along each $C_i$.
That is, there are divisors $D_i\subset X$
such that $D_i$ intersects $C_1+C_2$ only at $P_i\in C_i$ and
$D_i\cdot C_i=-K_X\cdot C_i$. Thus $D_1+D_2$ and $-K_X$ are
numerically equivalent in a neighborhood of $C_1+C_2$.
In particular, each $D_i$ intersects the general fiber of $f$ in a
single point and
$f:D_i\to S$ is an isomorphism near $0$. 

If $S$ is smooth at $0$ then $D_i$ is smooth at $P_i$.
We obtain a contradiction since then 
  $X$
is smooth at $P_i$. 
(If $0\in X$ is a 3--fold terminal singularity and $0\in H$ is a smooth
member of $|-K_X|$  then $X$ is smooth. This is quite easy to prove but
the only reference I know is the much too general 
\cite[3.5.1]{KoMo92}.)
Otherwise let
$\tilde S\to S$ be the universal cover of $S\setminus\{0\}$
and $\pi:\tilde X\to X$ the corresponding cover. 
$\pi$ ramifies only over $P_1$ and $P_2$, thus $\pi^{-1}(C_i)$
is a union of $\deg \pi$ curves intersecting at a single point over
$P_i$. Then $\pi^{-1}(C_1\cup C_2)$ is not a tree, a contradiction
again.\qed
\end{say}

\section{Irreducible fibers}

As (\ref{cart->cb.rem}) shows, we should  study the topology of conic
bundles and of rational curve fibrations with geometrically irreducible 
central fibers. In this section we prove that the geometrically
irreducible  cases are locally  quotients of conic bundles. This result
is   developed into a complete classification in the next section.

The proof relies on computing the algebraic fundamental group of
$X\setminus \sing X$ in  the neighborhood of a fiber.
First   I recall the definition and basic properties of
the algebraic fundamental group  of a scheme over $\r$.

\begin{say}[The algebraic fundamental group of a real   variety]
\label{real.fg.say}{\ }

As a general reference, see \cite{SGA1}.

For a complex variety $Y_{\c}$, 
the algebraic fundamental group
$\pi_1^{alg}(Y_{\c})$ is the profinite
completion of $\pi_1(Y(\c))$ (see, for instance, \cite[XII.5.2]{SGA1}).
If
$Y_{\r}$ is a  real variety, 
then $\pi_1^{alg}(Y_{\r})$  is related to $\pi_1(Y(\c))$
as follows. $\spec_{\r}\c\to \spec_{\r} \r$ is a degree 2 \'etale
morphism, and correspondingly we obtain a degree 2 \'etale cover
$Y_{\r}\times_{\r}\spec_{\r}\c\to \spec Y_{\r}$. 
The left hand side is essentially $Y_{\c}$, but viewed as  a real
variety. It has 2 geometrically connected components
(assuming that $Y$ itself is geometrically connected)
and the two components are interchanged by complex conjugation.
In particular, it has no real points.
Any further \'etale cover is obtained as  an \'etale cover of one of the
components, out of which we create the conjugate cover over the other
component.

Hence we get an exact sequence (cf.\ \cite[IX.6.1]{SGA1})
$$
1\to \hat{\pi}_1(Y(\c))\to \pi_1^{alg}(Y_{\r})\to \z_2\to 0.
$$
We see that if  an \'etale cover $\bar Y_{\c}\to Y_{\c}$
corresponds to a subgroup $H\subset \hat{\pi}_1(Y(\c))$, then the
choice of a real structure on $\bar Y_{\c}$ is equivalent
to a splitting of the above sequence modulo $H$. Such splittings need
not exist and they need not be unique.
\end{say}

We need the  computation of the real $\pi_1^{alg}$ in one case:

\begin{say}\label{et.cov.circle}[{\'Etale} covers of
$\r[x,y]/(x^2+y^2-1)$] {\ }

Adjoining $\sqrt{-1}$ we obtain $\c[x,y]/(x^2+y^2-1)$,
which is isomorphic to $\c^*$. $\pi^{alg}(\c^*)=\hat{\z}$, thus
the algebraic fundamental group of $\r[x,y]/(x^2+y^2-1)$ is an
extension of $\hat{\z}$ by $\z_2$. I claim that this extension is
trivial, that is
$$
\pi_1^{alg}(\r[x,y]/(x^2+y^2-1))\cong \hat{\z}+\z_2.
$$
To see this, we construct the tower of Galois extensions
corresponding to the finite quotients of $\hat{\z}$.

 Rotating by $2\pi/n$ acts on
 $\r[s_n,t_n]/(s_n^2+t_n^2-1)$ and  the ring of invariants
is generated by $x_n$ and $y_n$ as in (\ref{rot.invs.lem})
with a single relation $x_n^2+y_n^2=1$. 
Hence
$$
\r[x,y]/(x^2+y^2-1)\cong \r[x_n,y_n]/(x_n^2+y_n^2-1)
\subset \r[s_n,t_n]/(s_n^2+t_n^2-1)
$$
is a degree $n$ Galois extension with Galois group $\z_n$.
The corresponding map between the set of real points
is the $n$-sheeted cover $S^1\to S^1$.

If $n$ is odd then $\hat{\z}\onto \z_n$  and 
$\z_2\to \{1\}$ give the unique quotient 
of order $n$ of  $\hat{\z}+\z_2$. 
 For $n=2m$ even, there are two other quotients. One 
corresponds to
$\hat{\z}\onto \z_m$  and 
$\z_2\to \z_2$. This gives the extension
$$
\r[x_m,y_m]/(x_m^2+y_m^2-1)
\subset \c[s_m,t_m]/(s_m^2+t_m^2-1).
$$
Finally there is also a quotient $\hat{\z}\onto \z_{2m}$  and 
$\z_2\to \z_2\subset \z_{2m}$. This corresponds to the
Galois extension
$$
\r[x_n,y_n]/(x_n^2+y_n^2-1)
\subset \r[s_n,t_n]/(-s_n^2-t_n^2-1).
$$
Notice that  $\spec \r[s_n,t_n]/(-s_n^2-t_n^2-1)$ has no real points.

By contrast one should observe that the algebraic fundamental group
of $\r[x,x^{-1}]$ is the completion of the infinite dihedral group.
This is shown by the fact that the natural degree $n$ extension
$\r[t^n,t^{-n}]\subset \r[t, t^{-1}]$
is not Galois.
\end{say}

\begin{lem}\label{rot.invs.lem} Consider the $\z_n$-action on
$\r^2_{s,t}$  which is rotation by $2\pi/n$.
The ring of polynomial invariants is generated by 3 elements
$$
x_n:=\sum (-1)^j\binom{n}{2j}t^{2j}s^{n-2j},\quad
y_n:=\sum (-1)^j\binom{n}{2j+1}t^{2j+1}s^{n-2j-1}
$$
and $z:=s^2+t^2$, subject to the single relation
$$
x_n^2+y_n^2-z^n=0.
$$
\end{lem}

Proof. The action is
$$
s\mapsto \cos(2\pi/n)s+ \sin(2\pi/n)t,\quad
t\mapsto -\sin(2\pi/n)s+ \cos(2\pi/n)t.
$$
In order to understand the ring of invariants, note that under this
action
$$
s+\sqrt{-1}t\mapsto e^{2\pi\sqrt{-1}/n}(s+\sqrt{-1}t).
$$
Thus the invariants are $z:=s^2+t^2$ and the real and imaginary parts of
$(s+\sqrt{-1}t)^n$.  The latter are $x_n$ and $y_n$. 
Then $x_n^2+y_n^2=(s^2+t^2)^n=z^n$.\qed

\begin{say}[Trouble with terminology]{\ }

In describing the local structure of rational curve fibrations over
$\c$, it is natural to use the analytic topology on the base. Thus in
effect we
deal with triplets $f:X\to S$ where $S$ is a germ of a complex
analytic space, $X$ is a germ of a complex analytic space along a
compact set $f^{-1}(0)$  and
$f$ is a projective morphism. It is natural to do the same in the real
analytic case. The only problem is that the resulting objects
do not seem to have a standard   name.
(One could use the terminology ``complex analytic space with a real
structure".)

Let $Y_{\r}$ be a smooth real algebraic variety. Then $Y(\c)^{an}$ is a
complex analytic space with a real structure on it. This is not the same
as a real analytic space $Y(\r)^{an}$.  The real analytic space 
$Y(\r)^{an}$ corresponds only to the germ of $Y(\c)^{an}$ along
its real  points. It does not contain enough information to describe
$Y(\c)^{an}$ completely.

For us it is crucial to keep a neighborhood of the whole fiber, not
just a neighborhood of the real points of a fiber.
Thus if $f:X\to S$ is a rational curve fibration over $\r$
and $0\in S(\r)$ a  point then  first replace $S$ with a
small neighborhood of $0\in S^0\subset S(\c)^{an}$
and then replace $X$ by $X^0:=f^{-1}(S^0)\subset X(\c)^{an}$.

One can   think of these objects as    complex analytic
spaces with an antiholomorphic involution on them.
Thus the only problem is that they do not have a name.
\end{say}

\begin{notation}\label{cent.fib.not}
  Let $X$ be a real algebraic 3--fold with terminal
singularities such that $K_X$ is Cartier along $X(\r)$. 
Let $f:X\to S$ be a rational curve fibration and  $0\in S(\r)$  a point
such that $f^{-1}(0)$ is geometrically irreducible.

By (\ref{cart->cb.rem}) we know that either $K_X$ is Cartier along
$f^{-1}(0)$ or there is a unique conjugate pair of points of index $m>1$
along
$f^{-1}(0)$ (\ref{term.sings.defn}).

In the proof of the next theorem it is essential to distinguish
several structures on the central fiber $f^{-1}(0)$.
These are the following:
\begin{enumerate}
\item $f^{-1}(0)$, the central fiber as a real algebraic variety.
\item $\left(f^{-1}(0)\right)(\r)$, the topological space of real points
of the central fiber.
\item $f^{-1}(0)_{\c}$, the central fiber as a complex algebraic variety.
\item $\left(f^{-1}(0)\right)(\c)$, the complex analytic space of complex
points of the central fiber.
\end{enumerate}
There is not much problem in mixing up the last 2 objects,
but it is cruicial to keep the distinction between 
$f^{-1}(0)$ and $f^{-1}(0)_{\c}$ in mind.
\end{notation}

\begin{thm}\label{cover.thm}
 Notation and assumptions as in (\ref{cent.fib.not}).  There is a
real conic bundle $\tilde f:\tilde X\to \tilde S$ and an $\tilde
f$-equivariant
$\z_m$-action on $\tilde X$ and $\tilde S$ such that
$$
(f:X\to S)\qtq{is real analytically isomorphic to} (\tilde f:\tilde X\to
\tilde S)/\z_m.
$$
Moreover, we can choose $p_X:\tilde X\to X$ such that the following
conditions are satisfied:
\begin{enumerate}
\item  The induced morphism 
$$
S^1\sim \left(\tilde f^{-1}(0)\right)(\r)\to 
\left(f^{-1}(0)\right)(\r)\sim S^1
$$
is a degree $m$ cover.
\item If $s\in S(\r)\setminus\{0\}$, then $p^{-1}(f^{-1}(s))$ is
either isomorphic to $m$ copies of $f^{-1}(s)$
or $\left(p^{-1}(f^{-1}(s))\right)(\r)=\emptyset$ (the latter can happen
only for $m$ even).
\end{enumerate}
\end{thm}

Proof. Replace $S$ with a small analytic neighborhood of $0$. 
Let $P,\bar P\in \left(f^{-1}(0)\right)(\c)$ be the two points of
index $m$.
$$
\left(f^{-1}(0)\right)(\c)\setminus\{P,\bar P\}\cong \c^*\qtq{thus}
\pi_1(\left(f^{-1}(0)\right)(\c)\setminus\{P,\bar P\})\cong \z.
$$ 
By \cite[0.4.13.3]{Mori88},  the natural map
$$
\pi_1(\left(f^{-1}(0)\right)(\c)\setminus\{P,\bar P\})\to
\pi_1(X(\c)\setminus\{P,\bar P\})
$$
is surjective and its kernel is $m\z$. Thus 
$\pi_1(X(\c)\setminus\{P,\bar P\})\cong \z_m$.
Let $\tilde X_{\c}\setminus\{Q,\bar Q\}\to X_{\c}$ denote the universal
cover of 
$X_{\c}\setminus\{P,\bar P\}$; this can be extended uniquely to a
normal scheme
$\tilde X_{\c}$ which admits a finite morphism $\tilde X_{\c}\to
X_{\c}$.  $\tilde f_{\c}$ and $\tilde S_{\c}$ are obtained 
from the Stein factorization of $\tilde X_{\c}\to S_{\c}$.
$\tilde X_{\c}$ has only  index 1 terminal singularities
and so $\tilde f_{\c}$  is a conic bundle by (\ref{cart->cb.rem}).  All
that remains to do is to find a real algebraic variety $\tilde X$ such
that
$\tilde X_{\c}$ is the complexification of $\tilde  X$.

By (\ref{real.fg.say}), we have the following exact sequences, the
second obtained by taking quotient of the first by $m\hat{\z}$:
\begin{enumerate}
\setcounter{enumi}{2}
\item \quad
$0\to \hat{\z} \to \pi_1^{alg}(f^{-1}(0)\setminus\{P,\bar P\})
\to \z_2\to 0.$
\item \quad 
$0\to \z_m \to \pi_1^{alg}(X\setminus\{P,\bar P\})
\to \z_2\to 0.$
\end{enumerate}

By (\ref{et.cov.circle}) the above sequence (3) splits,
and so the same holds for (4). Moreover, in both cases we have a
distinguished splitting which induces a connected $m$-sheeted
covering of $\left( f^{-1}(0)\right)(\r)$.
The corresponding cover is denoted by $p_X:\tilde X\to X$.
This satisfies condition (\ref{cover.thm}.1) by construction.

Pick $s\in S(\r)\setminus\{0\}$. Then $\left(f^{-1}(s)\right)_{\c}$
is isomorphic to either $\c\p^1$ or    two copies of 
$\c\p^1$ meeting to a point. 
Thus   $\pi_1^{alg}(f^{-1}(s))=\z_2$.   $p^{-1}(f^{-1}(s))\to
f^{-1}(s)$ is \'etale and Galois. Thus $p^{-1}(f^{-1}(s))$
is either 
$m$ copies of
$f^{-1}(s)$ or $m/2$ copies of the nontrivial degree 2 cover
of $f^{-1}(s)$. The latter has no real points by
(\ref{real.fg.say}).\qed

\section{$\z_m$-actions on real conic bundles}

In this section  we complete the analysis of rational curve fibrations
with an irreducible central fiber. (\ref{cover.thm}) reduced this
question to the study of $\z_m$-actions on conic bundles.
We start by recalling the representation theory of $\z_m$ over $\r$,
mostly to fix notation.

\begin{notation}\label{rotation.not}
 The following are the irreducible real representations
of the cyclic group $\z_m$.
\begin{enumerate}
\item the 1--dimensional trivial representation ${\bold 1}$ (also denoted
by
${\bold 1}^+$),
\item if $m$ is even, we have the 1--dimensional sign representation
${\bold 1}^-$,
\item for any   $a\in \z$ we have a 2--dimensional 
representation $R_{a,m}$ (rotation by $a\frac{2\pi}{m}$), given as
$$
u\mapsto \cos(a\textstyle{\frac{2\pi}{m}})u+
\sin(a\textstyle{\frac{2\pi}{m}})v,\quad 
v\mapsto -\sin(a\textstyle{\frac{2\pi}{m}})u+
\cos(a\textstyle{\frac{2\pi}{m}})v.
$$
 $R_{a,m}\cong  R_{b,m}$  iff $a\equiv\pm b\mod m$.  $R_{a,m}$ is
irreducible iff
$a\not\equiv 0,m/2\mod m$. 
\end{enumerate}
${\bold 1}^{\pm}(z)$ denotes the vectorspace $\r  z$ with
$\z_m$-action ${\bold 1}^{\pm}$. 
$R_{a,m}(x,y)$ denotes that we have the above described rotation action
on the vector space $\r x+\r y$ in the basis $x,y$. 
\end{notation}

\begin{say}\label{rotation.say}
 We will need the 3--dimensional representations
of $\z_m$ modulo tensoring by ${\bold 1}^{\pm}$. 
The 3--dimensional representations are
\begin{enumerate}
\item ${\bold 1}^{\pm}\oplus {\bold 1}^{\pm}\oplus {\bold 1}^{\pm}$, and
\item ${\bold 1}^{\pm}\oplus R_{a,m}$ for $a\not\equiv 0,m/2\mod m$.
\end{enumerate}
By tensoring with  ${\bold 1}^-$ if necessary, these can be 
brought to the form
${\bold 1}(z)\oplus R_{a,m}(x,y)$ where $0\leq a<m$ is arbitrary
and $z,x,y$ is a choice of a basis. Under this action the
vector space of quadratic forms decomposes as follows:
$$
{\bold 1}(z^2)\oplus {\bold 1}(x^2+y^2) \oplus
R_{a,m}(zx,zy) \oplus
R_{2a,m}(x^2-y^2, 2xy).
$$
Let now $q(x,y,z)$ be a $\z_m$-equivariant real quadratic form
in the above variables
which is not the product of two linear forms over $\r$.
Assume that every nonidentity element of $\z_m$ has only isolated fixed
points on the conic $(q=0)$. Then $q$ is one of the following:
\begin{enumerate}
\setcounter{enumi}{2}
\item If $m\geq 3$ then  $q=\alpha z^2+\beta(x^2+y^2)$ and $(a,m)=1$.
\item If $m=2$ then 
$q=\alpha z^2+q'(x,y)$ and $a=1$. 
\end{enumerate}
\end{say}

We need the following lemma whose proof is left to the reader.

\begin{lem}\label{decomp.lem}
 Notation as above. Then
$$
\begin{array}{l}
R_{b,m}(s,t)\otimes S^2({\bold 1}(z)+R_{a,m}(x,y))
\cong\\
\qquad R_{b,m}(sz^2,tz^2) + R_{b,m}(s(x^2+y^2),t(x^2+y^2)) +\\
\qquad R_{a+b,m}(z(sx-ty),z(sy+tx)) + R_{a-b,m}(z(sx+ty),z(sy-tx)) +\\
\qquad R_{2a+b,m}(s(x^2-y^2)-2txy,t(x^2-y^2)+2sxy) +\\
\qquad R_{2a-b,m}(s(x^2-y^2)+2txy,t(x^2-y^2)-2sxy).\qed
\end{array}
$$
\end{lem}

\begin{notation}\label{quotient.thm.prenot}
 Let $f:Y\to S$ be a real conic bundle. Assume that
$\z_m$ acts $f$-equivariantly on  $Y$ and $S$. Let $0\in S(\r)$ be a
fixed point.  $f_*\o_Y(-K_{Y/S})$ is  a rank 3 real vector bundle
with a $\z_m$-action. In a neighborhood of $0\in S$ we can choose a
real analytic linearization of the actions, and we have the following
standard form:

$S$ is the germ of a 2--dimensional representation $(s,t)$ of $\z_m$.
$V$ is a 3--dimensional representation $(x,y,z)$ of $\z_m$
and $f_*\o_Y(-K_{Y/S})\cong V\otimes_{\r}\o_S$. Thus $Y$ has a
$\z_m$-equivariant embedding $Y\subset \p^2_{x,y,z}\times S$
and $g$ is the second projection.
\end{notation}

\begin{thm}\label{quotient.thm}
 Notation as above. Assume that $m\geq 2$, 
$\left(f^{-1}(0)\right)(\r)\neq \emptyset$ and
every $1\neq k\in \z_m$
acts with  only isolated and nonreal fixed points.
Then
$f:Y\to S$ is real analytically equivalent to one of the following
normal forms: 
\begin{enumerate}
\item  $m\geq 2$ and $Y=(z^2-x^2-y^2+h_1(x,y,z,s,t)=0)$, or
\item  $m$ is odd, $Y=(z^2+sx^2+2txy-sy^2+h_2(x,y,z,s,t)=0)$ and $b=2$.
\end{enumerate}
\noindent In both cases $V={\bold 1}(z)\oplus R_{1,m}(x,y)$, 
   $S=R_{b,m}(x,y)$ for some
$(b,m)=1$ and $h_i\in (s,t)^i$. 
\end{thm}

\begin{rem} It is not difficult to see that in both of the above cases,
after a further coordinate change
we can achieve that $h_i=0$.
Without the $\z_m$-action this is a standard result on conic bundles
(cf.\ \cite[Chap.I]{beau77}). 
\end{rem}

Proof. We are interested in the projectivization of the $\z_m$-action
on $f_*\o_Y(-K_{Y/S})\cong V\otimes_{\r}\o_S$, thus $V$ can be replaced
by $V\otimes {\bold 1}^-$ if necessary. Thus  we may assume that
$V={\bold 1}(z)\oplus R_{a,m}(x,y)$  for some $(a,m)=1$ by
(\ref{rotation.say}). 
By choosing the generator of $\z_m$ suitably, we can assume that $a=1$. 

Let $s\in S(\c)$ be a  fixed point of $k\in \z_m$. Then
$k$ has a fixed point on $f^{-1}(s)_{\c}$. Hence
the $\z_m$ action on $S$ has an isolated fixed
point at $0$.
Thus $0\in S$ is real analytically isomorphic to the germ of 
$R_{b,m}(s,t)$ for some
$(b,m)=1$.

Let $q(x,y,z)=0$ be the equation of the central fiber. Then
$q=\alpha z^2+\beta(x^2+y^2)$ if $m\geq 3$ and
$q=\alpha z^2+q'(x,y)$ if $m=2$ by (\ref{rotation.say}). 

Consider first the case when $(q=0)$ is a smooth conic with real points.
We can bring $q$ to the normal form 
$q=z^2-(x^2+y^2)$ if $m\geq 3$ and
$q= z^2-(x^2\pm y^2)$ if $m=2$. If $q= z^2-(x^2- y^2)$, then $\z_2$
has a real fixed point at $(1:1:0)$, which is not our situation.
This gives the case (\ref{quotient.thm}.1).

If $(q=0)$ is a pair of intersecting lines then their real intersection
point is fixed by $\z_m$, so this is impossible.

Finally consider the case when $(q=0)$ is a double line.
Then $q=z^2$, hence the equation of $Y$ can be written as
$$
z^2+sq_s(x,y,z)+tq_t(x,y,z)+(\mbox{higher order terms in $s,t$}),
$$
where $sq_s(x,y,z)+tq_t(x,y,z)$ is $\z_m$-invariant.
$Y$ has ony isolated singularities along the central fiber, thus
$z=q_s(x,y,z)=q_t(x,y,z)=0$ has only finitely many solutions.

From (\ref{decomp.lem}) we see that 
an expression of the form  $sq_s(x,y,z)+tq_t(x,y,z)$ can be
$\z_m$-invariant iff $b\equiv \pm a, \pm 2a\mod m$. In the former case
$q_s$ and $q_t$ are both divisible by $z$
(these are the $R_{a\pm b,m}$ summands in (\ref{decomp.lem})). 
Thus $b\equiv  \pm 2a\mod m$ and so $m$ is odd since $(b,m)=1$.
So $sq_s(x,y,z)+tq_t(x,y,z)$ is an element of one of the
summands $R_{2a\pm b,m}$ in (\ref{decomp.lem}). These can always be
brought to the form $sx^2+2txy-sy^2$ by a suitable change of the
$s,t$-coordinates.
\qed

\begin{lem}\label{one.comp.lem}
 Notation as in (\ref{quotient.thm}). Then
$(Y/\z_m)(\r)= (Y(\r))/\z_m$.
\end{lem}

This is of interest only if $m$ is even. By (\ref{rot.invs.lem}), the
quotient 
$S/\z_m$ is locally isomorphic to the singularity $(u^2+v^2-w^{m}=0)$,
which has 2 real components if $m$ is even. The main assertion of the
lemma is that our conic fibration has real points over only one of
these components.
\medskip

Proof. $\z_m$ acts freely on $Y(\r)$, so $Y/\z_m$  is smooth along 
$(Y(\r))/\z_m$. Thus $(Y(\r))/\z_m$ is a connected component of 
$(Y/\z_m)(\r)$.  We know that $(f^{-1}(0))(\r)\sim S^1$ is connected,
so $(Y/\z_m)(\r)$ is connected.\qed

\begin{cor}\label{over.one.comp}
Let  $X$ be a real projective 3-fold with terminal
singularities such that $K_X$ is Cartier along $X(\r)$.
Let    $f:X\to S$ be a    rational curve fibration over $\r$ such that
$-K_X$ is
$f$-ample. Let $N\subset X(\r)$ be a connected component.
Then $f(N)$ intersects only one of the connected components
of $S(\r)\setminus \sing S$.
\end{cor}

Proof. $f(N)$ is connected, thus if it 
intersects two connected components
of $S(\r)\setminus \sing S$ then there is a point
$P\in \sing S$ such that $P$ locally disconnects
$S(\r)$ and $f(N)$ intersects two of the local components of
$S(\r)$.  By (\ref{cart->cb.rem}), all real singular points of $S$ arise
as in (\ref{quotient.thm}). At each of these, even
$f(X(\r))$  intersects only one of the local components of 
$S(\r)\setminus \sing S$ by (\ref{one.comp.lem}), a contradiction.\qed

\section{Local topological description of real conic bundles}

In this section we study the topology of real conic bundles.
Aside from finitely many points the usual methods of conic bundle
theory (cf.\ \cite[Chap.I]{beau77}) give the following local analytic
description:

\begin{lem}\label{cb.top.gen.lem}
Let $S$ be a quasiprojective real algebraic surface
and $f:X\to S$ a real conic  bundle over $S$. 
There is a finite set of points $T\subset S(\r)$ such that every $s\in
S(\r)\setminus T$ has a Euclidean neighbourhood $s\in W_s\subset S$
such that $f:f^{-1}(W_s)\to W_s$ is fiber preserving real
analytically equivalent to one of the following normal forms.
(In all 4 cases $f$ is the
second projection and $B^2_{s,t}$ is the unit ball in $\a^2_{s,t}$.)
\begin{enumerate}
\item ($S^1$-bundle) 
$(x^2+y^2-z^2)\subset \p^2_{(x:y:z)}\times B^2_{s,t}$,
\item (empty fibers) 
$(x^2+y^2+z^2)\subset \p^2_{(x:y:z)}\times B^2_{s,t}$,
\item (collapsed end)
$(x^2+y^2+tz^2)\subset \p^2_{(x:y:z)}\times B^2_{s,t}$,
\item (blown up $S^1$-bundle)
$(x^2-y^2+tz^2)\subset \p^2_{(x:y:z)}\times B^2_{s,t}$.\qed
\end{enumerate}
\end{lem}

In order to understand $X(\r)$, we have to describe the local structure
of $f$ over the exceptional points $T$. This is quite difficult to do
up to real  analytic equivalence. It turns out, however, that  
if $X(\r)$ is a manifold, then a
suitable PL perturbation of $f$ becomes very simple everywhere.

\begin{thm}\label{loc.top.detailed.thm}
Let $S$ be a quasiprojective real algebraic surface
and $f:X\to S$ a real conic  bundle over $S$. 
Let $T\subset S(\r)$ be as in (\ref{cb.top.gen.lem}). 
Let $\bar n:\overline{X(\r)}\to X(\r)$ denote the topological
normalization.  Let  $M\subset \overline{X(\r)}$ be a connected component
which is a PL 3-manifold and set $U:=f(M)$. Let $T\cap U\subset
W_T\subset U$ be any open neighborhood.

Then there is a 2-manifold with
boundary $F\supset (U\setminus W_T)$ and a map $g:M\to F$
such that
\begin{enumerate}
\item $g|(M\setminus f^{-1}(W_T))=f|(M\setminus f^{-1}(W_T))$, and
\item   every $s\in F$ has a   neighbourhood $s\in W_s\subset F$
such that $g:g^{-1}(W_s)\to W_s$ is fiber preserving PL homeomorphic
  to one of the   normal forms (\ref{cb.top.gen.lem}.1--4). 
\end{enumerate}
\noindent Moreover, if $\bar n(M)$ is not a connected component of
$X(\r)$, then $g:M\to F$ contains  collapsed ends.
\end{thm}

Proof. 
  Let $n:\bar U\to U$ denote the normalization. As we noted
during the proof of (\ref{top.gen.cb.thm}), $f$ lifts to a morphism
between the normalizations $\bar f:M\to \bar U$. 
The locus of blown up fibers gives a curve $B\subset \bar U$.
We know that there is a finite set $\bar T:=n^{-1}(T\cap U)\subset \bar
U$ such that over $\bar U\setminus \bar T$ the map $\bar f$ is described
by one of the   normal forms (\ref{cb.top.gen.lem}.1--4). 

We study the various possibilities for the points $P\in \bar T$.

First consider the case
when $f^{-1}(n(P))$ is  reducible. Then   in suitable local coordinates
$X$ can be written as
$(x^2-y^2+h(s,t)z^2=0)$ for some real power series $h(s,t)$.
Along the fiber, $X$ has a unique singular point $Q$ at $x=y=0$
with local equation $(x^2-y^2+h(s,t)=0)$.

By \cite[4.4]{rat1}, the normalization of $X(\r)$ has a manifold
component near $Q$
iff one of the following conditions holds:
\begin{enumerate}
\setcounter{enumi}{2}
\item $(0,0)$ is an isolated point of $(h=0)$. Possibly after
interchanging $x$ and $y$, we may assume that $h(s,t)>0$ near $(0,0)$.
The projection
$$
\overline{X(\r)}\to \r\p^1_{(x:z)}\times \r^2_{s,t}
$$
is a 2-sheeted unramified cover. Thus $\overline{X(\r)}\to 
\r^2_{s,t}$  is an $S^1$-bundle.
\item $(h=0)$ is a PL 1--manifold near $(0,0)$. Then, up to a PL
coordinate change, we can assume that $h(s,t)$ is one of the local
coordinates. Thus, locally we have the standard equation 
 $(x^2-y^2+tz^2=0)$ for $X$.
\end{enumerate}

In all other cases, $\left(f^{-1}(n(P))\right)(\r)$ is $S^1$ or a point.
In particular, in all other cases, $\bar f^{-1}(P)$ is either $S^1$
or is collapsible to a point.
Let $P\in W_P\subset \bar U$ be a   regular neighborhood of $P$ such
that $N_P:=\bar f^{-1}(W_P)$ is a regular neighborhood of $\bar
f^{-1}(P)$. Note that $W_P$ is a disc if $P\in \inter \bar U$
and a half disc (with $P$ a boundary point) if $P\in \partial \bar U$.

If $P\in \inter \bar U$ then
$\partial N_P$ is an $S^1$-bundle over a circle with a point blown up for
each intersection point $B\cap \partial W_P$.
If $P\in \partial \bar U$ then
$\partial N$ is   $S^2$  with a point blown up for
each intersection point $B\cap \partial W_P$.

We consider separately the various topological types for $N_P$.

If $\bar f^{-1}(P)$ is 
 collapsible to a point, then $N_P$ is a ball, hence orientable.
Therefore,
$P\not\in B$ and $P\in \partial \bar U$. After collapsing
$\bar f^{-1}(P)$ to a point we have a collapsed end 
as in (\ref{cb.top.gen.lem}.3). 

Assume next that $\bar f^{-1}(P)\sim S^1$ and 
$N_P$ is a solid Klein bottle. If $P\in \inter \bar U$ then
 $P\not\in B$.  As explained in (\ref{def.tor-stor}.2), we can modify
$\bar f$ to get a map to an annulus with a collapsed end along the
inner circle. If $P\in \partial \bar U$ then
$\partial N_P$ is   $S^2$   blown up at 2 points. So $B\cap W_P$
is an interval and we can move $B$ away from
$P$. Thus we get disjoint curves of collapsed ends and blown up
$S^1$-bundles.

Finally consider the case when  $\bar f^{-1}(P)\sim S^1$ and 
$N_P$ is a solid  torus. Then $P\not\in \partial \bar U$
(for such a point   $\partial N_P$ would be $S^2$ or $S^2$ blown up
in a  few points) and  $P\not\in B$
(for such a point $\partial N_P$ would be a torus or Klein bottle
blown up in at least one point). 

Here we have to look at the algebraic structure of $f$. 
  If $f^{-1}(n(P))$ is a smooth conic, then $X(\r)\to S(\r)$
is an $S^1$-bundle near $n(P)$ and $M\to X(\r)$ is a homeomorphism near
$f^{-1}(n(P))$.

 We are left with the case when
$f^{-1}(n(P))$ is a double line.  Set $C=\red f^{-1}(n(P))$. Then
$C\cong \p^1$ and $(C\cdot K_X)=-1$. 
$S^1\sim \bar f^{-1}(P)\to C(\r)$ is a homeomorphism, thus
 $\overline{X(\r)}$ is not orientable along  $f^{-1}(P)$ by
(\ref{sing.nonorient.lem}), a contradiction.

The last assertion of the theorem is equivalent to the following local
statement which we checked in each case:

Let $n(P)\in D_P\subset S(\r)$ be a small neighborhood.
If $(f\circ \bar n)^{-1}(D_P)$ is disconnected then $M$
contains a collapsed end.
\qed

 \begin{lem}\label{sing.nonorient.lem}
 Let $X$ be a normal real algebraic variety and $C\subset X$
an irreducible  real algebraic curve which is not contained in  $\sing
X$. Assume that $K_X$ is Cartier along $C$ and $(C\cdot K_X)$ is odd.
Let $n:\overline{X(\r)}\to X(\r)$ denote the normalization.
Let $\gamma\subset  \overline{X(\r)}$ be a simple closed loop
such that $n:\gamma\to C(\r)$ is a homeomorphism.

Then $\overline{X(\r)}$ is not orientable along small perturbations of
$\gamma$.
\end{lem} 

Proof.  We can throw away the closed subset of $X$ along which $K_X$ is
not Cartier.  Let $\omega$ be a meromorphic real section of $K_X$
which is nonzero  and defined at all points of $C(\r)\cap \sing X$. Since
complex zeros and poles of $\omega|C$ come in conjugate pairs, we see
that the number of zeros of $\omega|C(\r)$ minus the number of poles of
$\omega|C(\r)$ is odd.  $n^*(\omega|X(\r))$ restricts to the
orientation bundle on the set of smooth points
(cf.\ \cite[2.7]{ras}).  Let $\gamma'$ be a  small perturbation of
$\gamma$ which is contained in the smooth locus. Then 
$n^*(\omega|X(\r))$ has an odd number of sign changes along $\gamma'$,
thus $\overline{X(\r)}$ is not orientable along $\gamma'$.\qed
\medskip

Here are some examples  of real conic bundles where the central
fiber is a double line.

\begin{exmp}  All of the examples sit in $\p^2\times \a^2$, $f$ is the
second projection and $U:=f(X(\r))$.

$(z^2+sx^2+2txy-sy^2=0)$.  $X$ is smooth and 
$f:X(\r)\to \r^2$ is an $S^1$-bundle outside
the origin.

$(z^2+sx^2+ty^2=0)$. $f:X(\r)\to \r^2$ 
has empty fibers over the positive quadrant and 
is an $S^1$-bundle over the other 3 open quadrants.  We have a
collapsed end along the positive $s$ and $t$ axes and a blown up
$S^1$-bundle along the negative $s$ and $t$ axes. 

$(z^2+sx^2+(t^2-s^2)y^2=0)$. $f:X(\r)\to \r^2$ 
has empty fibers over  $0<s<|t|$ and 
is an $S^1$-bundle outside
$st(s^2-t^2)=0$ otherwise.  
$U$ has 2 components: $(s\leq 0)$ and $(|t|\leq s)$.
We have a colapsed end everywhere along the $(|t|= s)$
curve. In the $(s\leq 0)$ component, there is a  blown up
$S^1$-bundle along the curve $B=(s+|t|=0)$. 
In the normalization map $\overline{X(\r)}\to X(\r)$ the preimage
of the central fiber is the disjoint union of a circle and of a point.

$(z^2+s^2x^2+2txy+s^2y^2=0)$.  $U$ has two components,
$(t\geq s^2)$ and $(-t\geq s^2)$. No blown up $S^1$-bundle curve. 
In the normalization map $\overline{X(\r)}\to X(\r)$ the central fiber
breaks up into 2 intervals.
\end{exmp}

\section{Proof of the main theorems}

The next theorem  contains all the information about
the topology of rational curve fibrations over $\r$
that we have obtained.  Later we
translate this information  to the language of the  classification
scheme  of 3--manifolds.

\begin{thm}\label{loc.top.detailed.gen.thm}
Let  $X$ be a real projective 3-fold with terminal
singularities such that $K_X$ is Cartier along $X(\r)$.
Let    $f:X\to S$ be a    rational curve fibration over $\r$ such that
$-K_X$ is
$f$-ample.
Let  $\bar n:\overline{X(\r)}\to X(\r)$ be the topological
normalization and $M\subset
\overline{X(\r)}$   a connected component which is a PL 3-manifold.
Then there is a 2-manifold with
boundary $F$ and a map $g:M\to F$
such that every $s\in F$ has a   neighbourhood $s\in W_s\subset F$
such that $g:g^{-1}(W_s)\to W_s$ is fiber preserving PL homeomorphic
  to one of the  following normal forms.
(In all 4 cases $g$ is the
second projection and $B^2_{s,t}$ is the unit ball in $\a^2_{s,t}$.)
\begin{enumerate}
\item ($S^1$-bundle) 
$(x^2+y^2-z^2)\subset \p^2_{(x:y:z)}\times B^2_{s,t}$,
\item (collapsed end)
$(x^2+y^2+tz^2)\subset \p^2_{(x:y:z)}\times B^2_{s,t}$,
\item (blown up $S^1$-bundle)
$(x^2-y^2+tz^2)\subset \p^2_{(x:y:z)}\times B^2_{s,t}$.
\item (Seifert fiber) 
$(x^2+y^2-z^2)/\z_m\subset (\p^2_{(x:y:z)}\times B^2_{s,t})/\z_m$
where $\z_m$ acts by rotation with angle $2\pi/m$ on $(x,y)$,
fixes $z$ and acts by rotation with angle $2b\pi/m$ on $(s,t)$
for some $(b,m)=1$.
\end{enumerate}
\noindent Moreover, the following additional conditions are also
satisfied: 
\begin{enumerate}
\setcounter{enumi}{4}
\item If $\bar n(M)$ is not a connected component of
$X(\r)$ then $g:M\to F$ contains  collapsed ends.
\item There is an injection   from the set
of  Seifert fibers of $g$ to  the set of singular points of $S$
contained in $f(M)$
with real analytic equation $(u^2+v^2-w^m=0)$. Under this injection the
multiplicity of the Seifert fiber equals the  exponent of $w$ in the
corresponding equation.
\end{enumerate}
\end{thm}

Proof. Let $P_j\in S(\r)$ be the set of points
such that $f$ is not a conic bundle over any  neighborhood of $P_j$.
This set is finite and (\ref{loc.top.detailed.thm}) describes the
behaviour of $f$ over $S\setminus\{\cup_jP_j\}$. 
Over suitable small neighborhoods of  each $P_j$, $f$ is described by
(\ref{quotient.thm}). 
If   (\ref{quotient.thm}.1) holds then we obtain a Seifert fiber
by (\ref{seifert.def}). If   (\ref{quotient.thm}.2) holds then
by (\ref{sing.nonorient.lem}) and 
(\ref{def.kb-skb}) we obtain a  curve of collapsed ends.\qed

\begin{say}[Proof of (\ref{orient.main.thm})]{\ }

Let us run the minimal model program for $f:X\to S$ over $\r$
\cite[Sec.3]{rat2}. 
At the end we obtain:
\begin{enumerate}
\item  A real algebraic threefold $X^*$ with only terminal singularities
such that $K_X$ is Cartier along $X^*(\r)$ 
and  $\overline{X^*(\r)}$ is  a 3--manifold (see  \cite[1.9,
1.11]{rat2}).
\item 
 A morphism
$f^{(*)}:X^*\to S^*$ to a real algebraic surface $S^*$
such that $-K_{X^*}$ is $f^{(*)}$-ample. (There is a birational
morphism $S^*\to S$ but   they need not be equal.)
\end{enumerate}

Moreover, by \cite[1.1]{rat2}, $X(\r)$ is obtained from 
$\overline{X^*(\r)}$ by the following operations:
\begin{description}
\item[\ ] throwing away all isolated points of $\overline{X^*(\r)}$,
\item[\ ] taking connected sums of connected components,
\item[\ ] taking connected sum with $S^1\times S^2$,
\item[\ ] taking connected sum with $\r\p^3$.
\end{description}

Let $M\subset X(\r)$ be a connected component and $M^*\subset
X^*(\r)$ its image. $M^*$  is connected, but its preimage 
$\overline{M^*}\subset \overline{X^*(\r)}$ may have several connected
components, say $M_1,\dots, M_s$.  $M_j$ minus finitely many points is
homeomorphic to an open subset of $M$, so each $M_j$ is orientable.

We can apply (\ref{loc.top.detailed.gen.thm}) to each $M_j$ 
to obtain maps  $g_j:M_j\to F_j$. $M_j$ is orientable,
thus we do not have blown up $S^1$-bundles, so every fiber of
$g_j$ is either a circle or a point.
By (\ref{top.sf-ls.thm}), each $M_j$ is either Seifert fibered
over $F_j$ or is a connected sum of lens spaces.
If $s\geq 2$ then each $M_j$ contains collapsed ends
by (\ref{loc.top.detailed.thm}), thus  one of the following holds:
\begin{enumerate}
\setcounter{enumi}{2}
\item $M\sim M_1\ \#\ a\r\p^3 \ \#\ b(S^1\times S^2)$
and $g_1:M_1\to F_1$ is Seifert fibered, or
\item $M\sim M_1\ \#\ \cdots\ \#\  M_s\ \#\ a\r\p^3 \ \#\ b(S^1\times
S^2)$ end each $M_i$ is a connected sum of lens spaces.
\end{enumerate}
By (\ref{intro.local.seif.fib}) there is an injection
from the set of the multiple Seifert fibers in case (3) (resp. the lens
space summands in case (4)) to the set of singular points of $S$ on
$f^{(*)}(M^*)$. 
By (\ref{over.one.comp}), $f^{(*)}(M^*)$ intersects
only one of the components $\bar U$ of $\overline{S^*(\r)}$. 

Assume first that  $f^{(*)}(M^*)\subsetneq \bar U$. Then $M\to F$
contains a curve of collapsed ends, so  we are in case (4). 
A singularity of type $(u^2+v^2-w^2=0)$ 
on $\bar U$ gives rise to a 
connected summand $L_{1,2}\cong \r\p^3$ by (\ref{top.sf-ls.compl}). 
  If $S_{\c}$
is rational then   there are at most 
6 singular points of type $(u^2+v^2-w^m=0)$ with $m\geq 3$ on $\bar U$
by (\ref{dvsurf.main.thm}), thus 
there are at most 6 lens spaces $L_{p,q}$ with $q\geq 3$.

Next consider the case when $f^{(*)}(M^*)= \bar U$. 
By (\ref{over.one.comp}) $\bar U$  can not map to
both local components of a singularity of type $(u^2+v^2-w^{2m}=0)$,
thus every such singular point of $\bar U$ is separating 
(\ref{separating.defn}). 
If $S_{\c}$
is rational then   there are at most 
6 singular points of type $(u^2+v^2-w^m=0)$ with $m\geq 2$ on $\bar U$
by (\ref{dvsurf.main.thm}), thus  $M\to F$ has at most 6 multiple
fibers.

An orientable component
of $\overline{S^*(\r)}$ is either  a sphere or a torus by
(\ref{dvsurf.main.thm}).\qed
\end{say}

In the nonorientable case the topological description  is more
complicated. The theorem below gives  substantial restrictions.
The proof yields a more precise result but I am not
sure how much further can one go.

\begin{thm}\label{nonorient.top.thm}
Let  $X$ be a real projective 3-fold with terminal
singularities such that $K_X$ is Cartier along $X(\r)$.
Let    $f:X\to S$ be a    rational curve fibration over $\r$ such that
$-K_X$ is
$f$-ample.
Let  $M\subset \overline{X(\r)}$ be a connected component
which is a PL 3-manifold. Then
$$
M\sim N\ \#\ L_1\ \#\ \cdots\ \#\  L_s\ \#\  a(S^1\times S^2)
\ \#\  b(S^1\tilde{\times} S^2),
$$
where each each $L_i$ is a lens space 
and  all the pieces of the 
Jaco--Johannson--Shalen decomposition of $N$ (cf.\ \cite[p.483]{Scott83})
are Seifert fibered.
\end{thm}

Proof.  By (\ref{loc.top.detailed.gen.thm}) there is a map
$g:M\to F$ whose local structure is given by
(\ref{loc.top.detailed.gen.thm}.1--4). 
In $F$ we have boundary curves $C_1,\dots,C_m$ and
curves corresponding to blown up $S^1$-bundles $B_1,\dots,B_k$.

We start with the boundary curves. As in the proof of
(\ref{top.sf-ls.thm}), we cut $F$ along simple paths starting and
ending in $\cup C_j$ which are disjoint from $\cup B_i$. Each time we
obtain either a connected sum decomposition or we remove a 1--handle.
At the end we may assume that each $C_j$ is parallel to a $B_i$. 
Let $B'_i$ be the boundary of a  regular neighborhood 
$U_i\supset B_i$. We assume that the $B'_i$ are disjoint. Let us cut $F$
along all the
$B'_i$. We obtain the pieces $U_i$ and some other pieces $V_k$.
$g^{-1}(V_k)\to V_k$ is Seifert fibered.
Each $U_i$ is either an annulus or a M\"obius band
so can be relaized as an interval bundle over a circle. Correspondingly,
 we obtain fibrations $g^{-1}(U_i)\to S^1$
whose fibers are annuli with one point blown up. Equivalently, these
fibers can be thought of as 
$$
R:=\r\p^2\setminus\{(1:0:0),(0:1:0)\}.
$$
Fibrations over $S^1$ with this fiber are classified by diffeomorphisms
of $R$ modulo isotopy. It turns out there there are only 8 such,
and they can be realized by the matrices
$$
\left(
\begin{array}{rr}
\pm 1 & 0\\
0 & \pm 1
\end{array}
\right)
\qtq{and}
\left(
\begin{array}{rr}
0 & \pm 1\\
\pm 1 & 0
\end{array}
\right)
$$
acting on the first 2 coordinates on $R\subset \r\p^2$.  It is easy to
see that each of these 8 cases gives  a Seifert fibered 3--manifold with
$\h^2\times
\r^1$-geometry. It is important to note   that
this Seifert fibration is not compatible with the Seifert fibration of
the pieces $g^{-1}(V_k)$.

The preimages of the $B'_i$ in $M$ are tori or Klein bottles,
so this starts to look like the Jaco--Johannson--Shalen decomposition.
We have to be careful, since we may have cut too much.

First, we can have a component $V_k$ which is a disc with  at most one
multiple Seifert fiber in it. Then   $g^{-1}(V_k)$ is a solid torus  
which can be pasted to the neighbouring $g^{-1}(U_j)$.
If $V_k$ does not contain any multiple  fibers then we get a collapsed
end and we can further decompose $M$ as above.
If $V_k$  contains a multiple  fiber  the we get a Seifert
fiber space. If there is no such $V_k$ then the cutting tori and Klein
bottles are all incompressible (cf.\ \cite[p.432]{Scott83}). 

Second, we can have a component $V_k$
which is an annulus with no multiple fibers. If one of the boundary
components of
$V_k$ is $B'_j$ the other $C_i$ then  over $V_k$ we
have a solid torus or Klein bottle which can be pasted to 
$g^{-1}(U_j)$ to obtain a new  Seifert fibered 3--manifold.
If the two boundary components of
$V_k$ are in $B'_j$ and $B'_{j'}$ then 
$g^{-1}(U_j), g^{-1}(U_{j'})$ and $g^{-1}(V_k)$ paste together
to a Seifert fibered 3--manifold with $\h^2\times
\r^1$-geometry. 

After all these operations we have Seifert fibered 3--manifolds
pasted along incompressible tori and Klein bottles.
Furthermore, no two pieces can be glued together to get a Seifert
fibered manifold (cf. \cite[p.440]{Scott83}). So this is the
Jaco--Johannson--Shalen decomposition.\qed

\section{Real surfaces with Du Val singularities}

The aim of this section is to 
study the configurations of Du Val singularities that can
simultaneously occur on a real algebraic surface which is rational over
$\c$.

\begin{defn}\label{separating.defn}
Let $S$ be a real algebraic surface and $P\in S(\r)$ a singular point
with equation $(u^2+v^2-w^{2m}=0)$. $\overline{S(\r)}$ has two
connected components locally near $P$. We say that $P$
is {\it separating} if these 2 local components are on different
connected components of $\overline{S(\r)}$ and {\it nonseparating}
if the 2 local components are on the same
connected component of  $\overline{S(\r)}$.
\end{defn}

\begin{thm}\label{dvsurf.main.thm}
 Let $S$ be a    projective real algebraic
surface such that $S_{\c}$ is rational. Assume that along $S(\r)$ we have
only singularities of real analytic type $(x^2+y^2-z^m=0)$. 
Let $M\subset \overline{S(\r)}$ be a connected component.
\begin{enumerate}
\item  If $M$ is orientable then $M\sim S^2$ or $M\sim S^1\times S^1$.
\item  $M$ can contain arbitrary  many points 
which map to a nonseparating singular point of type 
$(x^2+y^2-z^2=0)$ but 
 at most 6 points which map to other singular points of $S$.
\end{enumerate}
\end{thm}

The minimal model theory of real surfaces has been studied in detail
in the papers \cite{Comessatti14, Silhol89, ras}. It is not
difficult to generalize these results to the case when we allow
singularities of the form $(x^2+y^2-z^m=0)$. Most of the proofs run
very close to the arguments in \cite{ras}, thus I will only
explain the  points where some differences arise.

For the minimal model theory the natural
 setting is to allow all  Du Val singularities (\ref{dv.defn}).
(See \cite{Reid85} or \cite[4.2]{KoMo98} for the relevant background on
Du Val singularities over $\c$.)
This will show the surprising difference between
the $(x^2+y^2-z^m=0)$ and $(x^2-y^2-z^m=0)$ singularities.
(For $m=2$ the two singularities are isomorphic
which leads to some complications.)

\begin{defn}\label{dv.defn}
  A real Du Val singularity is a surface singularity
$(0\in S)\subset (0\in \a^3)$ which is real analytically equivalent to
one of the following normal forms:
\begin{description}
\item[$A_n^+$] $(x^2+y^2-z^{n+1}=0)$ for $n\geq 1$,
\item[$A_n^-$] $(x^2-y^2-z^{n+1}=0)$  for $n\geq 1$,
\item[$A_n^{++}$] $(x^2+y^2+z^{n+1}=0)$ for $n$ odd,
\item[$D_n^+$] $(x^2+y^2z+z^{n-1}=0)$ for $n\geq 4$,
\item[$D_n^-$] $(x^2+y^2z-z^{n-1}=0)$ for $n\geq 4$,
\item[$E_6^+$] $(x^2+y^3+z^4=0)$,
\item[$E_6^-$] $(x^2+y^3-z^4=0)$,
\item[$E_7$] $(x^2+y^3+yz^3=0)$,
\item[$E_8$] $(x^2+y^3+z^5=0)$.
\end{description}
$A_1^+\cong A_1^-$ but otherwise the singularities with different name
are not isomorphic. It is easy to see that these are all the real forms
of the complex Du Val singularities. (My setting and notation are not the
same as in
\cite[I.17.1]{AGV85},  but the above list is easy to obtain from the
results there.)
\end{defn}

\begin{defn} Let $0\in S$ be a smooth real point and $x,y$ local
coordinates at $0$.  The $(1,m)$-blow up of $(x,y)$ is the surface
$S'\subset S\times \p^1_{u:v}$ given by equation $ux-vy^m=0$.
For $m=1$ this is the ordinary blow up. The $(1,m)$-blow up has a
unique singular point of type $A_{m-1}^-$ at $(0,0,0,1)$. 

If $P,\bar P\in S$ are smooth and conjugate complex points,
then we can choose conjugate coordinate systems to do a $(1,m)$-blow up
at both points. The result is again a real algebraic surface with a 
conjugate pair of $A_{m-1}$-points (for nonreal points the signs in
the equations do not matter). 
\end{defn}

\begin{lem}\label{wbup.top.descr} Let $p:S'\to S$ be the $(1,m)$-blow
up of a smooth real point. 
If $m$ is even, then a small perturbation of $p$ gives a homeomorphism
$\overline{S'(\r)}\to S(\r)$. If $m$ is odd then
$S'(\r)\sim S(\r)\ \#\ \r\p^2$.

In both cases a $(1,m)$-blow
up creates a nonseparating singular point.
\end{lem}

Proof. If $m$ is odd then $(x,y)\mapsto (x,y^m)$ is a homeomorphism,
so the  $(1,m)$-blow up is  homeomorphic to the ordinary blow up.

Assume that $m=2n$ is even. Then  $S'$ has a singular point
$(ux-y^{2n}=0)$. In the normalization this splits into 2 parts,
and correspondingly the preimage of the exceptional $S^1$
becomes an interval. Thus $p:\overline{S'(\r)}\to S(\r)$
contracts an interval to a point, hence suitable small 
perturbations of $p$ are homeomorphisms.\qed
\medskip

The minimal model program for real surfaces is explained
in \cite{ras}, see also \cite{Silhol89, KoMo98}. The general theory
applies equally to surfaces with Du Val singularities.  For the
applications the key point is the following description of the extremal
contractions:

\begin{thm}\label{surf.dv.mmp}
 Let $F$ be a   projective  surface  over
$\r$  with Du Val singularities
and $R\subset \nec{F}$ a $K_F$-negative extremal ray. Then $R$ can
be contracted
$f:F\to F'$, and we obtain  one of the following cases:
\begin{enumerate}
\item[(B)] (Birational)   $F'$ is a  projective 
surface  over
$\r$  with Du Val singularities  and $\rho(F')=\rho(F)-1$.
 $F$ is the $(1,m)$-blow up of $F'$ at  a smooth point $P$ of $F'$ for
some
$m\geq 1$. We have two cases:
\begin{enumerate}
\item  $P\in F'(\r)$ is a real point, or
\item  $P$ is a pair of conjugate points.
\end{enumerate}

\item[(C)] (Conic bundle) $B:=F'$ is a smooth curve, $\rho(F)=2$ and
$F\to B$ is a conic bundle.

\item[(D)] (Del Pezzo surface) $F'$ is a point, $\rho(F)=1$, $-K_F$ is
ample and  $\rho(F_{\c})\leq 9$.
\end{enumerate}
\end{thm}

Proof. The proof can  be put together from the appropriate pieces
in
\cite{Cutkosky88, Morrison85, ras}. I just explain the main point: Why
does the   list  of birational contractions involve only $A^-$-type
singularities?

Let $E\subset F$ be the exceptional curve. Let $G\to F$ be the minimal
resolution of the singularities of $F$ along $E$ and $G\to F'$ the
composition. The exceptional divisor of $G\to F'$ consists of
$-2$-curves (coming from curves exceptional over $F$) and $(-1)$-curves 
(the irreducible components of the birational transform of $E$).
We can factor $G\to F'$ by repeatedly contracting  $(-1)$-curves. Two
$(-1)$-curves can never intersect since then we would get an exceptional
curve with selfintersection zero after contracting one of them.
By looking at the various cases we see that
  each connected component
of $\ex(G\to F')$ contains a unique $(-1)$-curve. In fact, we have one
of the following configuration of curves:
$$
\stackrel{-1}{\circ} - \stackrel{-2}{\circ} - \dots  -
\stackrel{-2}{\circ}
$$
Consider   curve germs $C_x$ and $C_y$ intersecting the $(-1)$-curve on
the left (resp. the $(-2)$-curve on the
right) transversally at a smooth point. After we contract everything,
$C_x$ and $C_y$ become a pair of transversally intersecting curves
$C'_x$ and $C'_y$ on $F'$.  Choose coordinates such that 
$C'_x=(x=0)$ and $C'_y=(y=0)$. We see that $F\to F'$ is the
$(1,m)$-blow up where $m-1$ is the number of $(-2)$-curves above.
\qed
\medskip

Since birational contractions  can eliminate only 
nonseparating $cA^-$-type
singularities,  we obtain:

\begin{cor}\label{dv.does.not.ch}
 Let $F$ be a    projective 
surface  over
$\r$  with Du Val singularities and $F\to F^*$ the result of the MMP
over $\r$. Then there is a one--to--one correspondence
between the two sets:
\begin{enumerate}
\item  Real Du Val singular points of $F$
which are either not of type $A^-$ or  are separating of type $A^-$.
\item  Real Du Val singular points of $F^*$
which are either not of type $A^-$ or are separating of type $A^-$.\qed
\end{enumerate}
\end{cor}

\begin{cor}\label{at.most.8.pts}
 Let $F$ be a    projective 
surface  over
$\r$  with Du Val singularities. Assume that $F_{\c}$ is rational.
Then every connected component of $\overline{F(\r)}$ contains at most
6   real Du Val singular points
which are either not of type $A^-$ or  are separating of type $A^-$.
\end{cor}

Proof.  By (\ref{dv.does.not.ch}) it is sufficient to consider the
cases when $F$ is either  a conic bundle or a Del Pezzo surface.

In the latter case  it is sufficient to consider degree 1 Del Pezzo
surfaces, since every other can be made into degree 1 by blowing up a
few points. Every degree 1 Del Pezzo surface is a double cover
of a quadric cone $Q\subset \p^3$, ramified along a curve $B$ 
not pasing through the vertex
which is
a complete intersection of $Q$ and of a cubic $C$. $B$ has the maximum
number of singular points $6$ when $C$ is the union of 3 planes.
(\cite{furu86} contains a partial list of singular Del Pezzo surfaces.)

If $f:F\to B$ is a conic bundle, then the number of singular points of
$F(\r)$ is not bounded, so we have to analyze the connected components
of  $\overline{F(\r)}$.

The local structure of  $f:F\to B$ is described by one of the
following equations. Here $s$ is a local coordinate on $B$
and we write $F$ as a hypersurface in $\p^2_{x:y:z}\times \a^1_{s}$. 
\begin{enumerate}
\item  $(x^2+y^2\pm z^2=0)$,
\item   $(x^2+y^2\pm s^mz^2=0)$,
\item  $(x^2+s(y^2+z^2)=0)$,
\item  $(x^2+syz=0)$,
\item  $(x^2+s y^2\pm s^mz^2=0)$,
\end{enumerate}
The first one describes an $S^1$-bundle or the empty set.
In case of the second equation, $F(\r)$ lies   entirely in $(\mp s\geq
0)$ if $m$ is odd and $\overline{F(\r)}$ decomposes into two 
connected components if $m$ is even and the sign is negative.
In the third case there are no real singular points and $F(\r)$ lies  
entirely in
$(\mp s\geq 0)$. In the fourth case $\overline{F(\r)}$ decomposes into
two  connected components, one in $yz-s\geq 0$ and one in $yz+s\geq 0$.
Each of the components passes through the two singular points
$(0,1,0,0)$ and $(0,0,1,0)$.

In the last case we have a $D$-type singular point. If we have the plus
sign then $F(\r)$ lies in $(s\leq 0)$. In the other case the local
structure is easy to work out by projecting to the $(x=0)$ line bundle.
We obtain that $\overline{F(\r)}$ decomposes into
two  connected components, one for $(s\geq 0)$ and one for $(s\leq 0)$.

Thus we see that each connected component of $\overline{F(\r)}$ pases
through
at most 4 singular points.\qed

\begin{say}[Proof of (\ref{dvsurf.main.thm})]{\ }

We can resolve the nonreal singular points of $S$ without changing  the
set of real points. Thus we may assume that $S$ has only Du Val
singularities. Let $f:S\to S^*$ be the result of the minimal model
program.

The statement about the number of singular points now follows from
(\ref{at.most.8.pts}).

Let $0\in S(\r)$ be a singular point with equation
$(x^2+y^2-z^m=0)$ and   $p:S'\to S$   the blow up of $0\in S$.
Then $S'(\r)\to S(\r)$ is a homeomorphism for $m\geq 3$. 
Thus the question of orientablility can be reduced to the case when all
singular points  have equation
$(x^2+y^2-z^2=0)$.

Assume that   $M\subset \overline{S(\r)}$ is orientable. 
$M$ gives a connected component of $M^*\subset \overline{S^*(\r)}$,
and 
$M^*$ is homeomorphic to $M$ by (\ref{wbup.top.descr}).

If $S^*$ is a conic bundle, then we obtain  from the proof of
(\ref{at.most.8.pts})
that every
connected component of $\overline{S^*(\r)}$ is $S^2, S^1\times S^1,
\r\p^2$ or a Klein bottle (these are the only surfaces that map to
$S^1$ such that every fiber is $S^1$, collapses to a point or is empty).

We are left with the case when $S^*$ is a Del Pezzo surface.
For every singularity $(x^2+y^2-z^2=0)$
I choose  the deformation $(x^2+y^2-z^2+\epsilon=0)$.
If these local deformations are simultaneously realized by a global
deformation   $S^*_{\epsilon}$ of $S^*$, then  it is
easy to see that $\overline{S^*(\r)}\sim S^*_{\epsilon}(\r)$. Thus 
(\ref{dvsurf.main.thm}.1) follows from the corresponding
result of Comessatti in the  smooth  case
(cf.\ \cite{Comessatti14, Silhol89, ras}).

For the existence  of $S^*_{\epsilon}$ it is sufficient to consider the
case when $S^*$ has degree 1 or 2. In the degree 2 case
$S^*$ is realized as a double cover of $\p^2$ given by an
equation $u^2=f_4(x,y,z)$. We look at deformations of the form
$u^2=f_4(x,y,z)+\epsilon g_4(x,y,z)$. We need to choose  $g_4$ to have
prescribed signs at the singular points of $S^*$.  Given at
most 6 points in $\p^2$ with no 5 on a line, one can always find  a
quartic $g_4$ which has prescribed values at the points. The degree 1
case is similar.
\qed
\end{say}

\section{Examples}

The aim of this section is to present examples of 3--manifolds
which can be realized as the real points of   rational curve
fibrations over rational surfaces. In many cases our examples are
unirational over $\r$, and they are always    rationally
connected (cf.\ \cite[IV.3]{koll96}).

\begin{say}[Circle bundles]\label{circle.bundles}{\ }

Let $F$ be a compact topological surface  (without boundary).
Circle bundles over $F$ are classified as follows (cf.\
\cite[p.434]{Scott83}). 

Let $p:M\to F$ be a circle bundle.
We can assume that it has
structure group $O(2,\r)$, since the homeomorphism group of the circle
retracts to $O(2,\r)$. (This result goes back to Poincar\'e, see
\cite[Sec.\ 4]{wood} for a proof.)
 Using the standard
2--dimensional representation of $O(2,\r)$, this induces an
$\r^2$-bundle    $E\to X(\r)$ with structure group $O(2,\r)$. That is,
in addition to the
vector space structure, each fiber carries an inner product, unique up to
a positive multiplicative constant. By a partition of unity argument we
can choose  a continuously varying inner product. Thus
$M$  can be thought of as the unit circles of an
$\r^2$-bundle $E\to F$ with an inner product.
The first
Stiefel--Whitney class  $w_1(E)$ of $E$ gives the first invariant.
The secondary invariant is  an element of
$H^1(F,R^1p_*\z)$.  The latter group is $\z_2$, except when $w_1(E)$ is
the orientation class of $F$, in which case it is $\z$.  (If $F$ is
orientable and the structure group is $SO(2,\r)$, then 
$E$ is naturally a $\c$-bundle and
this invariant coincides
with the first Chern class.)

Assume now that $F=S(\r)$. We would like to realize every 
topological circle bundle as a smooth conic bundle,
or better, as an algebraic circle bundle.
The above arguments suggest that $H^1(S(\r),\z_2)$ should be
``algebraic" for this to be possible. This can indeed be made precise,
but for us the main point is  the converse:

\begin{prop} Let $S$ be a smooth, projective real algebraic
surface such that $H^1(S(\r),\z_2)$ is generated by the cohomology
classes of algebraic curves. Then every topological circle bundle over
$S(\r)$ is fiber preserving homeomorphic to an algebraic circle bundle.
\end{prop}

Proof.  Let $M\to S(\r)$ be a circle bundle. We can assume that $M$ is 
the unit circle bundle of an
$\r^2$-bundle $E$ with an inner product.
The choice of inner products is equivalent to a section
$\sigma_0$ of  $S^2 E^*$.

By \cite[12.5.3]{BCR87} and
our assumption on $H^1(S(\r),\z_2)$, there is a strongly algebraic
vector bundle
$F\to S(\r)$ which is topologically isomorphic to $E$. 
(Strongly algebraic means that it is generated by its global section,
cf. \cite[12.1.6--7]{BCR87}.) $S^2 F^*$ is also strongly algebraic
by  \cite[12.1.8]{BCR87}, hence   $\sigma_0$ can be
approximated by algebraic sections  $\sigma_t$ by \cite[12.3.2]{BCR87}.
If $t$ is near $0$, then $\sigma_t$ defines an inner product
on each fiber and so we get an algebraic circle bundle.\qed
\medskip

If $S_{\c}$ is rational then $H^1(S(\r),\z_2)$ is generated by the 
cohomology classes of algebraic curves by \cite[p.67]{Silhol89}, thus
we obtain:

\begin{cor} Let $S$ be a smooth, projective real algebraic
surface such that $S_{\c}$ is rational. Then every topological circle
bundle over
$S(\r)$ is fiber preserving homeomorphic to an algebraic circle
bundle.\qed
\end{cor}

\end{say}

\begin{say}[Manifolds with spherical geometry]{\ }

The standard $n$-sphere  is given by equation
$$
S^n=(x_1^2+\cdots+x_{n+1}^2=1)\subset \r^{n+1}.
$$
Its group of automorphisms is $O(n+1,\r)$, which acts by real algebraic
automorphisms. Thus every quotient by a finite subgroup is a
unirational real algebraic variety.

In dimensions 3 every finite subgroup of $O(4,\r)$ acting freely on
$S^3$ is conjugate to a subgroup which leaves the Hopf fibration
invariant (cf.\ \cite[4.10]{Scott83}).
The Hopf fibration is easiest to write down as a map to the Riemann
sphere
$$
(x_1,x_2,x_3,x_4)\mapsto (x_1+\sqrt{-1}x_2)/(x_3+\sqrt{-1}x_4).
$$
Combined with the inverse of the stereographic projection
$$
s+\sqrt{-1}t\mapsto \left(\frac{2s}{1+s^2+t^2},
\frac{2t}{1+s^2+t^2},\frac{1-s^2-t^2}{1+s^2+t^2}\right)
$$
we get a real algebraic map
$$
p:(x_0^2+x_1^2+x_2^2+x_3^2=1)\to  (u_1^2+u_2^2+u_3^2=1).
$$
The coordinate functions of $p$ are easy to work out but they are
somewhat messy. The main point is that the fibers of $p$ are conics.
Indeed, from the first representation one sees that the fiber over
$s+\sqrt{-1}t$ is  given by equations
$$
x_3^2+x_4^2=(1+s^2+t^2)^{-1},\quad
x_1=sx_3-tx_4,\quad x_2=sx_4+tx_3.
$$
From this we see that  the indeterminacy locus of $p$ on the complex
projective quadric consists of the pair of conjugate lines
$$
(x_0=x_2-\sqrt{-1}x_1=x_4-\sqrt{-1}x_3),\ 
 (x_0=x_2+\sqrt{-1}x_1=x_4+\sqrt{-1}x_3).
$$
After blowing them up, $p$ becomes a $\p^1$-bundle over the 
projective quadric $(u_1^2+u_2^2+u_3^2=u_0^2)$.

Any finite subgroup as above acts on this conic bundle, so we see that
every 3-manifold with spherical geometry can be realized as a
unirational real rational curve fibration. 

It would be interesting to study the rationality question of these
quotients.
This should be of interest even over $\c$.
\end{say}

\begin{say}[3-manifolds with Euclidean geometry]{\ }
 Consider the following 4 finite cyclic subgroups of
$GL(2,\z)$:
$$
\begin{array}{rr}
\left\langle\left(
\begin{array}{rr}
0 & -1\\
1 & 1
\end{array}
\right)\right\rangle\ 
\cong \z_6, &
\left\langle\left(
\begin{array}{rr}
0 & -1\\
1 & 0
\end{array}
\right)\right\rangle\ 
\cong \z_4, \\
\left\langle\left(
\begin{array}{rr}
-1 & -1\\
1 & 0
\end{array}
\right)\right\rangle
\cong \z_3, &
\left\langle\left(
\begin{array}{rr}
-1 & 0\\
0 & -1
\end{array}
\right)\right\rangle
\cong \z_2.
\end{array}
$$
$GL(2,\z)$ acts on the torus $S^1\times S^1$ as described in
(\ref{def.tor-stor}). We can realize $S^1\times S^1$ as the real
algebraic group
$$
W:=\spec \r[x_1,x_2]/(x_1^2+x_2^2-1)\times
\spec \r[y_1,y_2]/(y_1^2+y_2^2-1),
$$
and the action of $GL(2,\z)$ is by algebraic automorphisms.
(For instance, the generator of $\z_4$ acts as
$(x_1,x_2, y_1,y_2)\mapsto (y_1,-y_2,x_1,x_2)$.)

Take $W\times \spec \r[z_1,z_2]/(z_1^2+z_2^2-1)$. Act on the $W$ factor
 by one of the
above groups  $\z_m$ and on the  second factor by 
rotation with angle $2\pi/m$. 

The quotient $X_m:=(W\times \spec \r[z_1,z_2]/(z_1^2+z_2^2-1))/\z_m$ is
a unirational real variety.  The resulting quotient   $M_m:=(S^1\times
S^1\times S^1)/\z_m$ is one of the components of
   $X_m(\r)$.
$M_m$  is Seifert fibered over $S^2$. In each of the above cases
we have 3 or 4 multiple fibers with the following multiplicities:
$$
\z_6: (6,3,2),\quad 
\z_4: (4,4,2),\quad 
\z_3: (3,3,3),\quad 
\z_2: (2,2,2,2).
$$

In the $\z_2$-case we can take a further quotient.
Let another copy of $\z_2$  act on $W$ by
$(x_1,x_2, y_1,y_2)\mapsto (-x_1,-x_2,y_1,-y_2)$. This is orientation
reversing on the set of real points.  The quotient 
$W(\r)/\z_2\times \z_2$ is homeomorphic to $\r\p^2$. 
On $\spec \r[z_1,z_2]/(z_1^2+z_2^2-1)$ one can act either as the
identity
or as $(z_1,z_2)\mapsto (z_1,-z_2)$.  In both cases the quotient 
$(S^1\times
S^1\times S^1)/\z_2\times \z_2$ is 
Seifert fibered over $\r\p^2$ with 2 fibers of multiplicity 2.
The quotient is nonorientable in the fist case and orientable in the
second case.

These 6 examples, together with the $S^1$-bundles found in
(\ref{circle.bundles}) exhaust all classical Seifert fiber spaces with
Euclidean geometry (cf.\ \cite[p.446]{Scott83}).

I have not tried to decide if these varieties are rational or not.
\end{say}

\begin{say}[Manifolds with Euclidean geometry]{\ }

More generally, let  $M$ be any compact manifold with Euclidean
geometry.
That is, $M$ is the quotient of Euclidean $n$-space by a group
$\Gamma$ of isometries. By the theorem of Bieberbach (cf.\
\cite[8.26]{raghu72})
$\Gamma$ contains $n$ linearly independent translations, so $M$ is a
quotient of the flat $n$-torus by a finite subgroup of $GL(n,\z)$.
If we represent the $n$-torus as the real algebraic variety
$$
T_n:=\prod_{i=1}^n \spec  \r[x_i,y_i]/(x_i^2+y_i^2-1),
$$
then $GL(n,\z)$ acts by algebraic automorphisms. So every quotient
by a finite group is a unirational real algebraic  variety.
\end{say}

\noindent University of Utah, Salt Lake City UT 84112 

\begin{verbatim}kollar@math.utah.edu\end{verbatim}

\end{document}